\documentclass[ijoc,nonblindrev]{informs3} 

\OneAndAHalfSpacedXII 

\usepackage{booktabs}
\usepackage{array}
\usepackage{multirow}
\usepackage{subfigure}
\usepackage{graphbox}
\usepackage[ruled,linesnumbered,noend]{algorithm2e}

\usepackage{natbib}
 \bibpunct[, ]{(}{)}{,}{a}{}{,}%
 \def\bibfont{\small}%
\TheoremsNumberedThrough     

\EquationsNumberedThrough    

\begin{document}


\RUNAUTHOR{Wang and Hong}

\RUNTITLE{Large-Scale Inventory Optimization: A Recurrent-Neural-Networks-Inspired Simulation Approach}

\TITLE{Large-Scale Inventory Optimization: A Recurrent-Neural-Networks-Inspired Simulation Approach}

\ARTICLEAUTHORS{%
\AUTHOR{Tan Wang}
\AFF{School of Data Science, Fudan University, Shanghai 200433, China\\\EMAIL{dwang19@fudan.edu.cn}}
\AUTHOR{L. Jeff Hong\footnote{Corresponding author}}
\AFF{School of Management and School of Data Science, Fudan University, Shanghai 200433, China\\\EMAIL{hong\_liu@fudan.edu.cn}}

} 

\ABSTRACT{%
Many large-scale production networks include thousands types of final products and tens to hundreds thousands types of raw materials and intermediate products. These networks face complicated inventory management decisions, which are often too complicated for inventory models and too large for simulation models. In this paper, by combing efficient computational tools of recurrent neural networks (RNN) and the structural information of production networks, we propose a RNN inspired simulation approach that may be thousands times faster than existing simulation approach and is capable of solving large-scale inventory optimization problems in a reasonable amount of time. 
}%

\HISTORY{}

\maketitle

%

\section{Introduction}\label{sec:intro}

Inventory is one of the most important tools for mitigating uncertainty in enterprise production and operations. Effective management of inventory can offer tremendous potential for increasing manufacturing efficiency and reducing operational cost. Therefore, inventory management has always been a core part of modern enterprise supply chain management and has long captured the interests from academics as well as industrial practitioners.

\begin{figure}[ht]
\begin{center}
\includegraphics[scale=0.3]{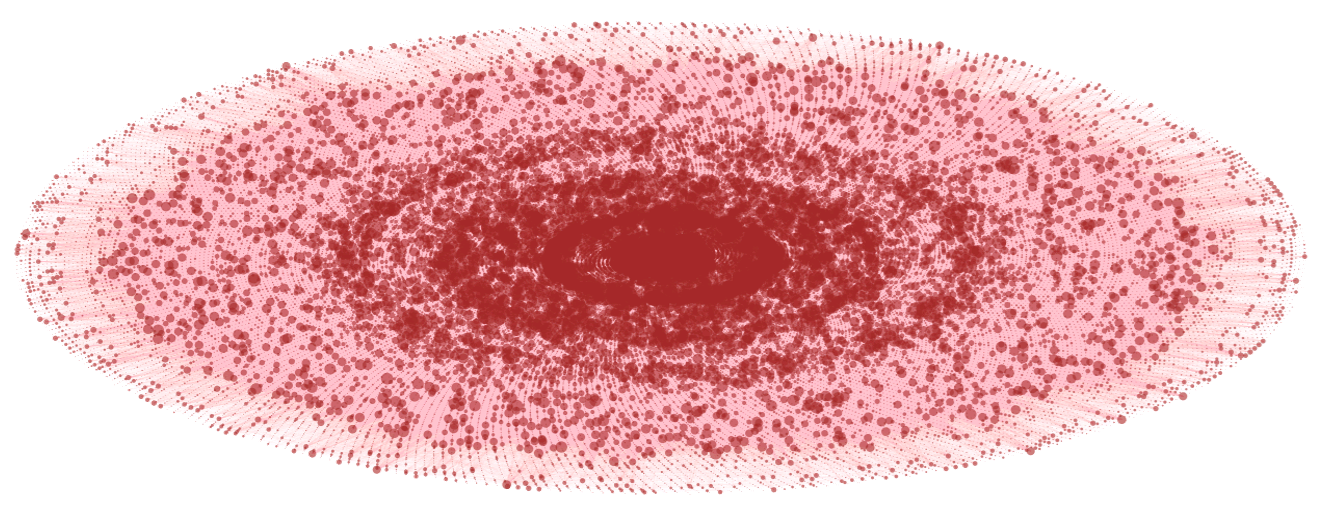}
\caption{Large Scale Inventory System} 
\label{fig:large scale system}
\end{center}
\end{figure}

Nowadays, many large companies produce and offer a wide range of brands or products, e.g., Coco-Cola and Unilever in fast-moving consumer goods industry and Samsung and Xiaomi in consumer electronics industry, partly due to market differentiation and partly due to economies of scale. As a result, the bill of materials (BOM) of these companies often include tens of thousands even to hundreds of thousands of nodes, which represent raw materials, sub-assemblies and final products, and show complex network topologies, which represent the production relationship of all the nodes. These complex production systems also impose tremendous challenges and opportunities to inventory management. First, there are decisions that need to be made for every node of the BOM, such as whether to keep inventory and how much to keep. Second, the complex network topology means a lot of sharing and pairing among the nodes in the production process and, to achieve the overall efficiency, one needs to take a holistic view towards all inventory decisions. Third, the overall inventory costs of these production systems are often very high and, therefore, any savings may be quite significant. For instance, in our consulting experience that motivated this study, we worked with a global manufacturer whose BOM has over $500,000$ nodes and over $4,000,000$ links (a part of the BOM is shown on Figure \ref{fig:large scale system})  and whose inventory at the time was worth about \$3 billion, thus even a $1\% $ decrease is significant. With the ongoing global pandemics and international trade frictions, the global supply chains are facing increasingly more risks and, therefore, how to manage these large-scale production and inventory systems becomes not only challenging but also extremely important.

From the perspective of research, the aforementioned problem belongs to multi-echelon (or multi-stage) inventory optimization problem. In the existing literature, two types of models have been developed to handle the problem: stochastic-service (SS) models and guaranteed-service (GS) models. The two types of models differ in replenishment mechanism between stages. SS models assume that the delivery or service time can vary based on the material availability at the supply stage, while GS models assume that each stage can quote a delivery or service time that can be always satisfied \citep{graves2003supply}. With regard to SS models, the existing work mainly focuses on serial systems \citep{clark1960optimal} or assembly systems \citep{rosling1989optimal,chen2014integrality}, which is in general not suitable for the large-scale problems considered in this paper. 

GS models were first proposed by \cite{simpson1958process}. \cite{graves2000optimizing} develop a dynamic programming algorithm for supply chains that can be modeled as spanning trees. \cite{lesnaia2005complexity}  show that the GS models for general acyclic networks are NP-hard problems and difficult to solve efficiently, and \cite{humair2011optimizing} propose two faster heuristics to solve the problems approximately. However, GS models assume that it is possible to establish meaningful deterministic upper bounds on the stochastic demands, which may be quite difficult for the large-scale problems considered in this paper and may lead to either uncontrolled fill rates or high inventory costs. One advantage of GS models is that they tend to hold inventory only at a small number of strategic nodes of the BOM. This property is particularly appealing for large-scale problems because it drastically reduces the difficulty of managing inventories, compared to the situation where all nodes hold inventory. In this paper, even though we do not use GS models, we want to develop algorithms that have this property.

The simulation approach is another approach to solving inventory problems, and it allows stochastic demands and complex BOM structures. Its general framework was laid out by \cite{glasserman1995sensitivity}. Basically, the simulation approach first builds a simulation model that simulates the evolution of the production and inventory system with randomly generated demands and given inventory policies, and estimates the average cost of the system by running the simulation model for multiple replications. To minimize the inventory cost, the simulation approach computes the sample-path gradient with respect to the inventory decisions and applies the stochastic approximation (SA) algorithm (also known as stochastic gradient descent (SGD) algorithm) to optimize.

Even though the simulation approach is simple conceptually and it is capable of modeling complex production-inventory systems, it has two drawbacks that prevent it from solving the large-scale problems considered in this paper. {\it First and foremost, it is often too slow for large-scale problems.} For instance, for a BOM network with 5000 nodes, we observe that the traditional simulation approach takes about an hour just to run a single simulation replication with 100 periods (see Table 1 in Section 6.1) and many hours to compute a sample-path gradient. Furthermore, we show that the computational complexities of the simulation and the gradient calculation are of $O(Tn^2)$ and $O(Tn^3)$, respectively, where $T$ is the number of periods and $n$ is the number of nodes in the BOM. Therefore, it is unlikely that the simulation approach may be used to solve problems with tens of thousands to hundreds of thousands nodes. {\it Second, it is not clear how to use the simulation approach to find solutions that keep inventory only at a small number of nodes.} By applying the SA algorithm directly, we observe that the solution typically keeps inventory for almost all nodes and it is very difficult to implement in practice.

To overcome the first drawback of the traditional simulation approach, we  notice that the simulation of a production-inventory system is very similar to the forward pass of a recurrent neural network (RNN), which is widely used in speech recognition and language modeling, and the simulation optimization of the inventory policy is similar to the training of a RNN. This analogy is critical to the speedup of the simulation approach, because RNNs are sometimes bigger than the production-inventory systems that we consider and the computational tools that are used to model and to train RNNs, e.g., the back-propagation (BP) algorithm for computing sample-path gradient and the various SGD algorithms for optimization, may be adapted to the simulation and optimization of the inventory decisions. Furthermore, we discover that, different from RNNs, BOM networks are typically quite sparse, which allows us to use sparse matrices to further speed up the simulation and optimization. We show that, by combining the tools of RNNs and the use of sparse matrices, we improve the computational complexities of simulation and gradient computation from $O(Tn^2)$ and $O(Tn^3)$ to $O(Tn)$ and $O(Tn)$ respectively, and reduce the computational times in the order of thousands to tens of thousands.

To overcome the second drawback of the traditional simulation approach, we propose to use $L_1$-regularization to force the solution of the optimization problem to be sparse, thus only allowing a small number of nodes to have non-zero base-stock levels. While regularization is commonly used in training of RNNs to avoid overfitting, it is used in our problems for completely different purposes. However, the similarity between the training of RNNs and the inventory optimization allows us to take full advantage of the sub-gradient SGD algorithm, designed for training $L_1$-regularized deep neural networks, to solve inventory optimization problems. Furthermore, we show that fast iterative shrinkage-thresholding (FISTA) algorithm, which utilizes the Nesterov accelerated gradient method and the proximal method to achieve faster convergence, may be applied to our problem as well and achieve numerical performance that is in general slightly better than the sub-gradient SGD algorithm. We also show that, because of the difference between the training of RNNs and the inventory optimization, a re-optimization step may be added to our problem to further improve the performance of the inventory decisions while keeping the inventory only at a small number of nodes.

After overcoming the two drawbacks of the traditional simulation approach, the new simulation approach is capable of solving large-scale inventory optimization problems with tens of thousands to hundreds of thousands nodes in a reasonable amount of time, while keeping the inventory only at a small number of nodes. Furthermore, we test different methods to implement the algorithms, including the use of TensorFlow for simulation and for automatic calculation of the sample-path gradients, the use of parallel CPU processors and the use of GPUs, and obtain a wide range of insights that are useful in implementing the algorithms in different scenarios.

In summary this paper makes the following contributions to large-scale inventory optimization as well as simulation modeling and optimization:
\begin{itemize}
\item It finds that the modeling of a complex production-inventory system is analogous to a RNN, and its simulation optimization is analogous to the training of a RNN. This analogy opens the door so that we can take advantages of efficient computational tools of deep learning to solve large-scale inventory problems. In addition, we want to point out that this analogy goes beyond production-inventory systems. It may be applied to general periodically reviewed dynamic systems and, thus, may lead to speedup of the simulation of general dynamic systems.

\item It demonstrates that the BP algorithm provides the same sample-path gradient as the infinitesimal perturbation analysis (IPA). However, it is computationally more efficient than the IPA when the gradient is of a high dimension. It also explores the possibility of using TensorFlow or other modern computational tools to compute the sample-path gradient automatically without deriving it explicitly.

\item It shows that the sparsity of the BOM network may be utilized to significantly reduce the computational complexities of the simulation and gradient calculation, allowing the algorithm to reduce the computation time by orders of magnitude.

\item It uses $L_1$-regularization to keep inventory only at a small number of nodes, thus reducing the difficulty of inventory management, and proposes algorithms to efficiently solve the $L_1$-regularized simulation optimization problems.
\end{itemize}

The rest of this paper is organized as follows. In Section \ref{sec:state}, we provide the problem statement and simulation-optimization formulation of the problem. In Section \ref{sec:rnn}, we briefly introduce the typical structure of RNNs and make the analogy between the modeling and training of RNNs and the simulation and optimization of large-scale inventory systems. Inspired by RNNs, in Sections \ref{sec:model} and \ref{sec:opt}, we discuss how to utilize the structures of inventory systems to develop efficient simulation and optimization algorithms; followed by a comprehensive numerical study in Section \ref{sec:exp}. The paper is concluded in Section \ref{sec:conclud}.

\section{Problem Formulation}\label{sec:state}

We consider a large-scale production system with a general BOM network structure, where the demands and inventories are reviewed periodically (e.g., every day). We differentiate the items of the BOM into two categories, procurement items and production items. Procurement items are the raw materials that sit at the bottom of the BOM, and they are procured from outside. Production items are the intermediate items (a.k.a. sub-assemblies) or final products, and they are produced by the production system. We allow all production items to have outside demands, because some intermediate items are needed for maintenance or service. In our model, the outside demands, the production time and the procurement times may be stochastic. But we assume that they follow known distributions (or, at least, may be generated through a simulation algorithm) that take only integer values.

Furthermore, we assume that there are no capacity constraints on the production. This is a common assumption in the inventory management literature \citep{clark1960optimal,rosling1989optimal,graves2000optimizing}. For large-scale production-inventory systems, considering capacity constraints will turn the inventory optimization problem into a complex production planning/scheduling problem \citep{hall2010capacity}, which requires separate optimization tools and typically only considers deterministic/known demands. In fact one may consider the uncapacitated inventory optimization problem studied in this paper as an add-on to the production plans to hedge the randomness in the demands.

In addition, we suppose that (installation) base-stock policies are used for all nodes of the BOM, and our goal is to find the optimal base-stock levels. Notice that the base-stock policies may not be the optimal policies for our problem. However, they are widely used in practice because of its simplicity \citep{glasserman1994stability,gallego1999stock}, and they are known to be optimal for serial systems \citep{clark1960optimal} and assembly systems \citep{rosling1989optimal}.

\subsection{The Simulation-based Approach}\label{subsec:simu-approach}

As the BOM network is complex and there are external demands for possibly all items in the BOM, the simulation algorithm also includes a lot of details. To clearly present the algorithm, we break it into the following four major steps based on the sequence of events within each period.
\begin{enumerate}
  \item[\bf Step 1.] 
  Calculating the inventory positions, observing outside demands and placing the replenishment orders at the beginning of the period.
  \item[\bf Step 2.]
  Receiving finished orders, fulfilling outside demands, and updating the on-hand inventories and backlogs at the beginning of the period. 
  \item[\bf Step 3.]
  Producing based on the replenishment orders and the available materials, and updating the on-hand inventories at the end of the period.
  \item[\bf Step 4.]
  Calculating the cost at the end of the period.
\end{enumerate}
Now we present the algorithm based on these four steps.

\subsubsection{Step 1.}
Let ${IP}_{t,i}$ denote the inventory position of item $i$ at period $t$. It is a measure used under a base-stock policy for making ordering decisions, which equals the on-hand inventory plus the on-order quantity (i.e., the amount that is ordered but not yet received through either procurement or production) minus the backorders (i.e., the demands that have occurred but have not been fulfilled) \citep{snyder2011fundamentals}. Let $O_{t,i}$, $D_{t,i}^{in}$ and $D_{t,i}^{out}$ denote the order quantity, the internal demand and the outside demand of item $i$ at period $t$, respectively. We show in the Appendix that  the inventory position can be calculated using the following recursive formula:
\begin{equation}
I{P_{t,i}} = I{P_{t - 1,i}} + {O_{t - 1,i}} - D_{t - 1,i}^{in} - D_{t - 1,i}^{out},
\label{eq.IPti-recursive}
\end{equation}
where $I{P_{0,i}}$ is the initial inventory position of item $i$ at the beginning of time $0$.

Once the inventory position $I{P_{t,i}}$ is calculated, a replenishment order $O_{t,i}$ is placed to bring $I{P_{t,i}} - D_{t,i}^{out} - D_{t,i}^{in}$ up to the base-stock level $S_i$, i.e.,
\[{O_{t,i}} =  - \min \left\{ {0,I{P_{t,i}} - D_{t,i}^{out} - D_{t,i}^{in} - {S_i}} \right\}.\]
Notice that, while the outside demands $D_{t,i}^{out}$ are clearly observed, the internal demands $D_{t,i}^{in}$ depend on the replenishment orders of downstream items. Therefore, the replenishment order of item $i$ can not be placed until orders of all its downstream items are placed. Let $a_{ij}$ denote the number of units of the component $i$ are required to produce a unit of item $j$. Then, the internal demand of item $i$ by its downstream item $j$ observed at period $t$ is ${a_{ij}}{O_{t,j}}$, and the total internal demand is $D_{t,i}^{in} = \sum\limits_j {{a_{ij}}{O_{t,j}}}$. Furthermore, these calculations must start from the final products, which are the most downstream items and have no internal demands, and gradually move to the upstream items.

\subsubsection{Step 2.}
Let $P_{t,i}$ denote the arrived replenishment orders of item $i$ at the beginning of period $t$. They may be procured from external suppliers if the item is a raw material item, or they may be produced within the production system if the item is a production item (i.e., a sub-assembly or a final product). Let $I_{t,i}$ and $B_{t,i}^{out}$ denote the on-hand inventory and backlogged outside demand at the {\it end} of period $t$, respectively. 

We assume the outside demands have higher priority and are fulfilled at the beginning of the period when the items are available. Then, the on-hand inventory at the beginning of period $t$ after fulfilling the outside demand, denoted by $I_{t,i}^0$, is updated by
\begin{equation}
I_{t,i}^0 = \max \left\{ {0,{I_{t - 1,i}} + {P_{t,i}} - B_{t - 1,i}^{out} - D_{t,i}^{out}} \right\},
\label{eq.Iti-start}
\end{equation}
and the unfilled part is backlogged with
\begin{equation}
B_{t,i}^{out} =  - \min \left\{ {0,{I_{t - 1,i}} + {P_{t,i}} - B_{t - 1,i}^{out} - D_{t,i}^{out}} \right\}.
\label{eq.Btiout}
\end{equation}

Let $l_{t,i}$ denote the (random) lead time of item $i$ ordered at time $t$. It may be the procurement lead time if the item is procured from outside, or the production lead time if the item is produced within the production system. Notice that, for a raw material item (say $i$), the procurement arrived at period $t+l_{t,i}$ is the order at period $t$, i.e.,
${P_{t + {l_{t,i}},i}} = {O_{t,i}}$. The production items are much more complicated, because it is determined not only by the lead time but also the availability of upstream items. We handle them in {\bf Step 3}.

\subsubsection{Step 3.}

For a production item, orders are manufacturing commands and they are processed immediately if the upstream components are all available. Notice that we do not assume that the BOM has a tree structure as commonly assumed in the literature \citep{graves2000optimizing, rosling1989optimal}. Then, some components may be shared by multiple downstream items. When the inventories of these components are insufficient, an allocation rule is needed. In this case, for simplicity, we assume that the inventory is allocated proportionally to all downstream demands so they all have the same fill rate. Notice that for any item, the shortage is a small probability event. Then, the use of different allocation rules may not have significant impact on the performance of the system. The proportional allocation rule is simple and easy to implement, it avoids solving complex allocation optimization problems, and it maintains sample-path continuity that is crucial to the calculation of sample-path gradient. Let $r_{t,i}$ denote the fill rate of item $i$ at period $t$. Then,
\[{r_{t,i}} = \min \left\{ {\frac{{I_{t,i}^0}}{{\sum\limits_j {{a_{ij}}\left( {{O_{t,j}} + O_{t - 1,j}^b} \right)} }},1} \right\},\]
where $j$ denotes an immediate downstream item of item $i$, $a_{ij}$ denotes the number of units of the component $i$ are required to produce a unit of item $j$ and $O_{t-1,j}^b$ is the backlogged production orders of item $j$ at period $t-1$ due to the unavailability of at least one of its upstream items.

Let $M_{t,i}$ denote the production quantity of item $i$ at period $t$. With the allocated upstream inventories and newly assigned production orders, we can now calculate the production quantity $M_{t,i}$ and the backlogged production order $O_{t,i}^b$ as follow:
\begin{eqnarray}
{M_{t,i}} &=& \mathop {\min }\limits_{j,{a_{i,j}} \ne 0} \left\{ {{r_{t,j}}} \right\}\left( {{O_{t,i}} + O_{t - 1,i}^b} \right),\label{eq.Mti}\\
O_{t,i}^b &=& {O_{t,i}} + O_{t - 1,i}^b - {M_{t,i}},\label{eq.Otib}
\end{eqnarray}
where the minimum operator used in Equation (\ref{eq.Mti}) takes into consideration of the lowest availability of all components that are necessary in producing item $i$. Notice that
${M_{t,i}}$ will become finished goods (or arrived order) after a production lead time of $l_{t,i}$ periods, i.e., 
\[{P_{t + {l_{t,i}},i}} = {M_{t,i}}.\]

Furthermore, for an item $i$, the total internal demand by its downstream item $j$ at period $t$ is ${a_{ij}}\left( {{O_{t,j}} + O_{t - 1,j}^b} \right)$. As some production orders may be backlogged due to insufficient supply, the actual amount to be deducted from the on-hand inventory of item $i$ is ${a_{ij}}{M_{t,j}}$. Therefore, the on-hand inventory of each item $i$ at the end of each period $t$ is
\begin{equation}
{I_{t,i}} = I_{t,i}^0 - \sum\limits_j {{a_{ij}}{M_{t,j}}}.
\label{eq.Iti-end}
\end{equation}

\subsubsection{Step 4.}

The cost of each period is composed of the inventory holding cost and the stockout penalty cost, i.e.,
\[{C_t} = \sum\limits_i {\left( {{h_i}{I_{t,i}} + {p_i}B_{t,i}^{out}} \right)}, \]
where ${{h_i}}$ and ${{p_i}}$ denote the unit holding cost and unit stockout penalty cost of item $i$, respectively. Notice here  we only consider the stockout penalty cost for backlogged outside demands.

\vspace{11pt}

\begin{figure}[ht]
\begin{center}
\includegraphics[scale=0.38]{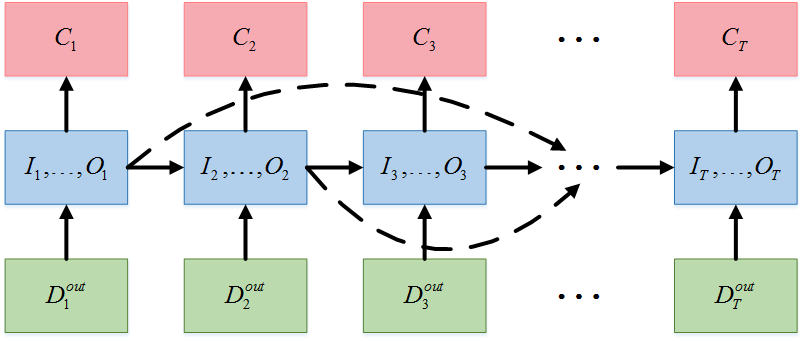}
\caption{Schematic Diagram of Inventory Simulation} 
\label{fig:simulation sketch}
\end{center}
\end{figure}

The aforementioned {\bf Steps 1} to {\bf 4} form an iteration (i.e., a period) of the inventory simulation algorithm, and it runs for $T$ periods to calculate (an observation of) the cumulative cost. We illustrate the algorithm in Figure \ref{fig:simulation sketch}.

The objective of the simulation-based approach is to find the optimal base-stock level that minimizes the expected cumulative cost over $T$ periods, i.e.,
\begin{equation}\label{eq:optprob}
\min\quad E\left[ {\sum_{t = 1}^T {{C_t}\left( S \right)} } \right].
\end{equation}
Here we write the cost of a period $t$ as $C_t(S)$ just to emphasize that the cost is a random function of the base-stock levels $S=(S_1,\ldots,S_n)$.

Because the simulation model itself is very complicated and there exists no approach to deriving a closed-form expression of the expected cumulative cost, we suggest to solve Problem (\ref{eq:optprob}) using the SA approach, which was also adopted by \cite{glasserman1995sensitivity} to solve a much smaller scale inventory optimization problem. The key to use the SA approach is to calculate the sample-path gradient of {\small${\sum\limits_{t = 1}^T {{C_t}\left( S \right)} }$} with respect to $S$, which can be solved by using the IPA and will be discussed in details in Section \ref{subsubsec:ipa}.

\subsection{Challenges in Large Scale Problems}\label{subsec:challange}

\begin{figure}[ht]
\begin{center}
\subfigure[Serial System]{ 
\raisebox{4.5\height}{
\includegraphics[align=c, width=0.3\columnwidth]{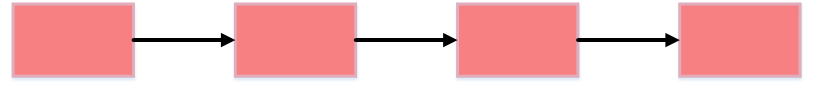}}
\label{fig:topology:a}
} 
\subfigure[Assembly System]{ 
\includegraphics[width=0.3\columnwidth]{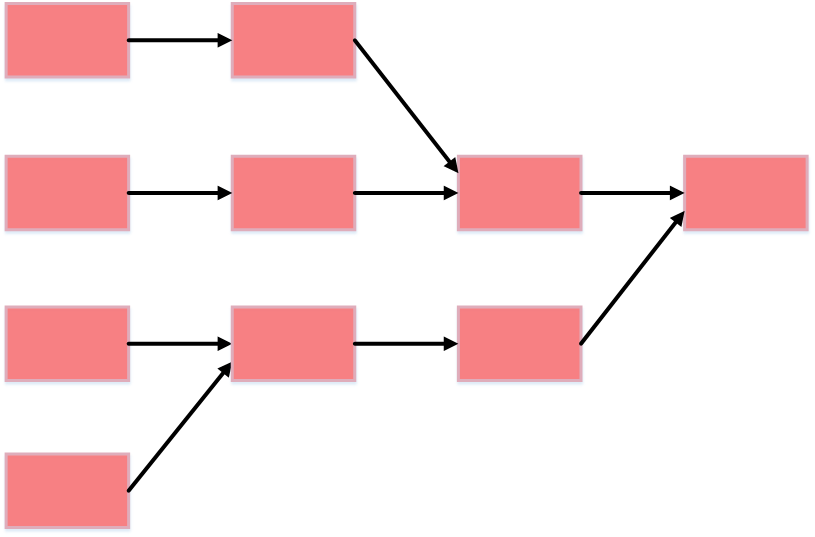}
\label{fig:topology:b}
}
\subfigure[General System]{ 
\includegraphics[width=0.3\columnwidth]{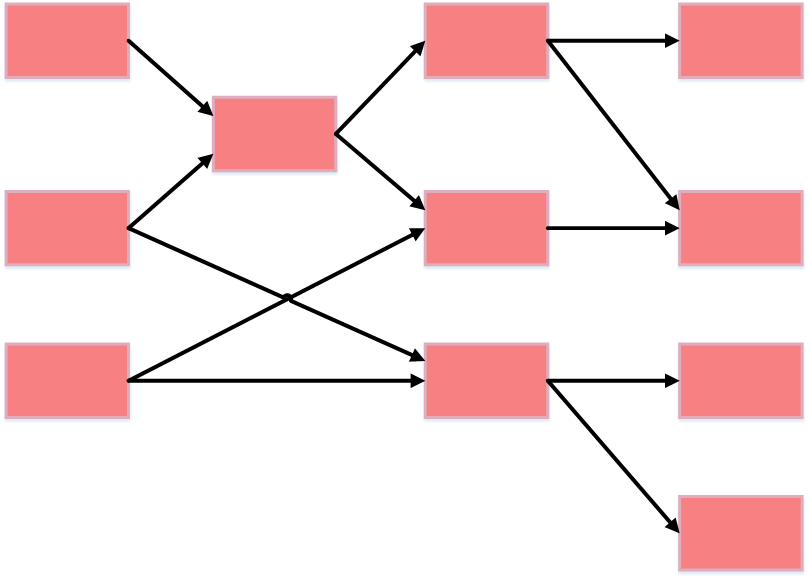}
\label{fig:topology:c}
}
\caption{Network Topologies \citep{snyder2011fundamentals}} 
\end{center}
\end{figure}

Previous studies on inventory simulation mainly focus on serial systems (see Figure \ref{fig:topology:a}) or assembly systems (see Figure \ref{fig:topology:b}), the simulation-based approach introduced in Section \ref{subsec:simu-approach} is capable of handling inventory systems with general topology (as illustrated in Figure \ref{fig:topology:c}). However, this generality also brings tremendous challenges in computation, especially for large-scale problems that are considered in this paper.

Notice that in serial systems, each stage has at most one predecessor and at most one successor, and in assembly systems, each stage has at most one successor. Therefore, the computation needed to simulate and optimize these systems does not increase dramatically as the size of the problem increases. In general systems, however, the complex network structure implies that, for every node of the BOM, any other node may be its predecessors or successor. This structure significantly complicates the simulation process and results in a drastically increase of the simulation time when the problem is of large scale. 

For instance, according to our testing experience, the simulation of a single replication of a system with only $5,000$ nodes for $100$ periods may take about an hour, and the calculation of its sample-path gradient may take over ten hours. Moreover, while multiple simulation replications are easy to parallel, the parallelization of a single replication of the simulation run or the parallelization of a single gradient calculation is not straight-forward and may need careful redesign of the simulation algorithm. Therefore, the first challenge of this paper is {\it how to reduce the computational time (including both simulation and optimization) to a manageable level for large-scale problems.}

Furthermore, when the system is large and with tens of thousands or even more items, keeping inventory for all items may not be practical. However, directly applying the SA algorithms to solve Problem (\ref{eq:optprob}) may result in solutions that are nonzero for almost all items and, thus, very difficult to implement in practice. One may add a cardinality constraint on the number of nonzero base-stock levels, but it is very difficult to handle in the context of simulation optimization. Therefore, the second challenge of this paper is {\it how to find solutions that keeps inventory only at a small number of nodes.}

\section{Recurrent Neural Networks}\label{sec:rnn}

RNN is a class of neural networks that allow previous outputs to be used as inputs while having hidden states. This enables it to exhibit temporal dynamic behaviors and to model time series data. RNN models have achieved the-state-of-the-art performance on tasks that include speech recognition \citep{miao2015eesen}, language modeling \citep{mikolov2010recurrent}, and image captioning \citep{you2016image}. The scale of deep RNN models can be very large, for example, the model of \cite{miao2015eesen} contains 8.5 million parameters.

\subsection{RNN Architecture and Training}\label{subsec:rnn-archit}

\begin{figure}[ht]
\begin{center}
\includegraphics[scale=0.38]{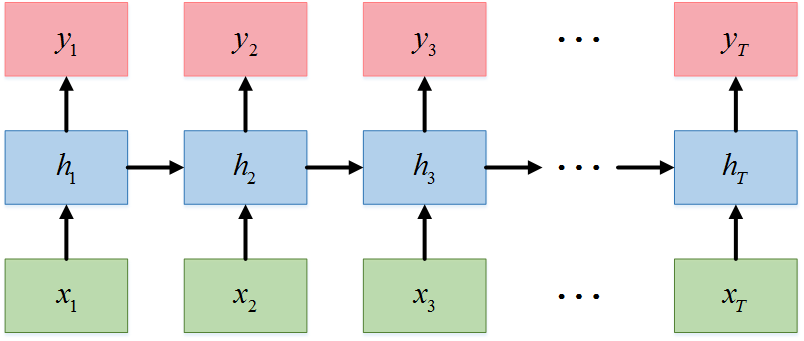}
\caption{Typical RNN Architecture} 
\label{fig:rnn structure}
\end{center}
\end{figure}

As shown in Figure \ref{fig:rnn structure}, an RNN often has a chain-like architecture, where in each period, the middle box is a neural network that takes $x_t$ and $h_{t-1}$ as the inputs and outputs $h_t$ and $y_t$. While many variations exist, a common implementation of an RNN can be described by the following updating equations \citep{goodfellow2016deep}:
\begin{eqnarray}
{h_t} &=& f\left( {W{h_{t - 1}} + U{x_t} + b} \right), \label{eq:rnn1}\\
{y_t} &=& g\left( {V{h_t} + c} \right), \label{eq:rnn2}
\end{eqnarray}
where $h_t$ denotes the hidden states at time $t$, $f(\cdot)$, $g(\cdot)$ are complex nonlinear activation functions represented by neural networks and $U,V,W,b,c$ are the parameters of the neural networks.

Suppose that we have a dataset of $N$ samples, denoted by $(x_t^n,y_t^n)$, $t=1,2,\ldots,T$ and $n=1,2,\ldots,N$. Given the network parameters, with an initial value of the hidden state $h_1$ and the outside inputs $\{x_1^n,x_2^n,\ldots,x_T^n\}$, we can use the RNN model to generate a time series of outputs, denoted by $\{\hat y_1^n,\hat y_2^n,\ldots,\hat y_T^n\}$. Notice that, if the neural network perfectly matches the data generating process, then $\hat y_t^n = y_t^n$. We define a loss function $L\left( {y_t^n,\hat y_t^n}  \right)$ to capture the difference between $\hat y_t^n$ and $y_t^n$. Then, to train an RNN is to find the network parameters that solves the following optimization problem:
\begin{equation}\label{eq:rnntraining}
\min\quad \frac{1}{N}\sum\limits_{n = 1}^N {L_T^n}, \quad {\rm where}\quad L_T^n = \sum_{t = 1}^T L\left( {y_t^n,\hat y_t^n}  \right).
\end{equation}

To solve Problem (\ref{eq:rnntraining}), SGD algorithms are typically used. Notice that SGD algorithms are originated from SA algorithms, and they have been adapted to solve very high-dimensional problems \citep{bottou2018optimization}. To use SGD algorithms, one needs to calculate the gradient of $L_T^n$ with respect to the network parameters. This is typically done through BP algorithm. The BP algorithm for calculating the gradient of neural networks is a celebrated result in the history of deep learning \citep{rumelhart1986learning}. It is particularly efficient if the network parameters are of high dimension, and it may be implemented automatically \citep{baydin2018automatic}. The particular version of the BP algorithm for RNNs is called back-propagation through time \citep{werbos1990backpropagation}, and it is capable of handling neural networks with a few million parameters \citep{miao2015eesen}.

\subsection{Similarities Between RNN and Inventory Optimization}\label{subsec:rnn-similarity}

There are three critical similarities between the training of RNN and the simulation approach introduced in Section \ref{sec:state}. First, both models have a complex network that evolves over time, and both have external and internal inputs in each period. In the case of RNNs, the network is a neural network, the external and internal inputs are $x_t$ and $h_{t - 1}$ respectively, as illustrated in Figure \ref{fig:rnn structure}. In the case of inventory model, the network is the production network implied by the BOM, the external and internal inputs are $D_t^{out}$ and $I_{t - 1},\ldots,O_{t - 1}$ respectively, as illustrated in Figure \ref{fig:simulation sketch}. Because the lead times may not always be one period, the inventory model may have links that go beyond the immediate next period, as illustrated in Figure \ref{fig:simulation sketch}. Therefore, the simulation model is more complex than the RNN model in this aspect.

Second, both problems have very similar optimization formulations. Notice that the objective function of the inventory optimization problem is an expectation, which may be approximated by the sample average of $N$ samples. Then, the problem becomes
\[\min\quad \frac{1}{N}\sum\limits_{n = 1}^N {\sum\limits_{t = 1}^T {C_t^n\left( S \right)} },\]

which has the exact same form as the RNN training problem (\ref{eq:rnntraining}). Furthermore, the decision variables of the inventory optimization problem are the base-stock levels. They can also be considered as the parameters of the production network. Therefore, the decision variables of both problems are network parameters that do not change over time.

Third, both approaches calculate sample-path gradient and optimize using gradient-based algorithms. The training of RNN uses the BP algorithm for gradient calculation and the SGD algorithm for optimization; while the inventory optimization problem uses the IPA algorithm for gradient calculation and the SA algorithm for optimization.

Given all these similarities between the two problems, it is natural to ask {\it ``how the training of RNN is capable of solving problems with millions of decision variables?''} and {\it ``what can be learned to solve the large-scale inventory optimization problems?''}. These are the questions that we answer in next sections.

\section{Simulation Model}\label{sec:model}

RNN and other neural network models are typically represented in terms of tensors. The concept of a tensor is a generalization of vectors and matrices and can be easily understood as a multidimensional array. Many computational tools, such as TensorFlow and PyTorch, have been developed for fast tensor operations, and they can easily leverage on the parallel computing capability provided by multi-core CPUs and many-core GPUs to speed up the calculations. Indeed, tensors are so important that Google even named its machine learning library as ``TensorFlow". Therefore, to achieve a significant speedup of the simulation model for large-scale inventory problems, a critical step that we learned from the success of RNN is the tensorization of the simulation model.

\subsection{Matrix Representation of BOM}\label{subsec:bom-net}

A BOM is a list of raw materials, sub-assemblies and their quantities and manufacturing relations needed to manufacture the final products. Because we consider general inventory systems in this paper, the BOMs of these systems may be represented by complex directed networks. Therefore, we adopt the adjacency matrix that is commonly used in network science and graph theory to represent the BOM. The adjacency matrix of a BOM with $n$ items is a square matrix with $n\times n$ elements, denoted by $A$, where its $(i,j)^{th}$ element $a_{ij}$ indicates that $a_{ij}$ units of item $i$ are required to produce a unit of (downstream) item $j$. From the viewpoint of a network, it indicates that there is a directed arc from item $i$ to item $j$ and the weight of the arc is $a_{ij}$ if $a_{ij}>0$. If there is no arc, then $a_{ij}=0$.

Notice that the adjacency matrix has $n \times n$ elements. For large-scale supply chains with $10^4$ to $10^5$ items, the storage of the matrices can take up a lot of memory space and the matrix operations also need significant amount of computational power. Nevertheless, according to our observations of practical large-scale supply chains, the adjacency matrix is typically very sparse. In fact, sparsity is a common feature of large-scale complex networks in the real world \citep{ugander2011anatomy,wang2012network}. Therefore, we use sparse matrices to represent the adjacency matrices to save memory space as well as to improve the computation efficiency.

For a directed network with $n$ nodes, its density is defined by
\[\rho  = \frac{m}{{n\left( {n - 1} \right)}},\]
where $m$ is the actual number of edges \citep{albert2002statistical}. Notice that if the network is fully connected, the actual number of edges is $n(n-1)$. Therefore, the density is a measure of the network sparsity. Given the density $\rho$, the average (in) degree of the network, denoted by $\left\langle k \right\rangle $, is 
\[\left\langle k \right\rangle  = \frac{m}{n} = \left( {n - 1} \right)\rho  \approx n\rho .\]
In the context of a BOM network, the average (in) degree indicates the average number of upstream components needed to assemble a downstream product, which typically does not grow with the scale of the network. Therefore, it is reasonable to assume that the $\left\langle k \right\rangle $ is a constant, i.e., 
\begin{equation}\label{eqn:sparsity}
n\rho\sim {\rm constant}. 
\end{equation}
This relationship is important for our understanding of the impact of the sparsity on the performance of the simulation algorithm as well as the inventory optimization algorithms.

\subsection{Tensorization of the Simulation Model}\label{subsec:model-imple}

Notice that in the simulation model presented in Section \ref{subsec:simu-approach}, except for the adjacency matrix, all other variables at time $t$ may be represented by vectors with $n$ elements, e.g., the on-hand inventory vector ${I_t} = \left( {{I_{t,1}}, \ldots ,{I_{t,n}}} \right)$ and the order quantity vector ${O_t} = \left( {{O_{t,1}}, \ldots ,{O_{t,n}}} \right)$ and so on. With these vectors and the adjacency matrix, the equations introduced in Section \ref{subsec:simu-approach} can be rewritten in the form of vector and matrix operations. Indeed, most of these transformations are straightforward, here we only present a few operations related to the adjacency matrix.

\begin{itemize}
\item {\bf Replenishment Orders.} Notice that the internal demand of an item depends on the orders of its downstream items. It may be calculated by $D_t^{in} = {O_t} \times {A^T}$. Therefore, the update of inventory positions (Equation \ref{eq.IPti-recursive}) may be written as
\[I{P_t} = I{P_{t - 1}} + {O_{t - 1}} - {O_{t - 1}} \times {A^T} - D_{t - 1}^{out}.\]

Furthermore, the replenishment order of an item depends on the internal demands that depend on the replenishment orders of the the downstream items. Therefore, in Section \ref{subsec:simu-approach}, we proposed an iterative algorithm that starts from the most downstream items (i.e., final products) to calculate the replenishment orders of all items. Notice that the number of iterations depend on the number of layers of the BOM network. Suppose that there are ${n_l}$ layers. Then, the longest distance between two nodes (items) is ${{n_l}-1}$. We know that for a directed network with an adjacency matrix denoted by $A$, ${\left( {{A^r}} \right)_{ij}} > 0$ if and only if it takes $r$ steps from node $i$ to node $j$ \citep{wang2012network}. Here, ${{A^r}}$ denotes the r-th power of matrix $A$. Thus, the number of layers can be calculated by 
\[{n_l} = \sup \left\{ r \in N, \sum\limits_{i,j} 1_{ \left\{\left( A^r \right)_{ij} > 0\right\}}  > 0 \right\} + 1,\]
where $1_{\{\cdot\}}$ is an indicator function. 

\item {\bf Fill Rates.} The fill rate vector of the items can be calculated by 
\[ r_t = \min\left\{I_t^0 / \left[ \left (O_t+O_{t-1}^b \right) \times A^T \right], 1 \right\}
\]

where the operators ``/'' and $\min (  \cdot  )$ denote the element-wise division and minimization. \end{itemize}

\vspace{17pt}

Let $C_t^{sum} = \sum\limits_{u = 1}^t {{C_u}} $ denote cumulative cost up to period $t$ and let $I^{initial}$ denote the given initial on-hand inventory. The tensor representation of the simulation model is presented in Algorithm \ref{alg:simulation}. It is worth noting that as the individual production quantity depends on the lowest availability of all necessary components and the lead time of each item is different. Therefore, there are still some operations that cannot be tensorized in the algorithm.

\begin{algorithm}
\small
\label{alg:simulation}
\SetAlgoNoLine
\caption{The Tensorized Simulation Algorithm}
\KwSty{Initialization:} $I{P_0} = {I_0} = {I^{initial}},{O_0} = 0,D_0^{out} = 0,B_0^{out} = 0,O_0^b = 0,C_0^{sum} = 0$\;
\For {$t = 1$ \KwTo $T$}{
$I{P_t} \leftarrow I{P_{t - 1}} + {O_{t - 1}} - {O_{t - 1}} \times {A^T} - D_{t - 1}^{out}$ \label{line-simu:IPt}\;
${O_t} \leftarrow  - \min \left\{ {0,I{P_t} - D_t^{out} - S} \right\}$ \label{line-simu:Ot1}\;
\For {$l = 1$ \KwTo ${n_l-1}$}{
$IP_t^{temp} \leftarrow I{P_t} - D_t^{out} - {O_t} \times {A^T}$ \label{line-simu:IPttemp}\;
${O_t} \leftarrow  - \min \left\{ {0,IP_t^{temp} - S} \right\}$ \label{line-simu:Ot2}\;
}
$I_t^0 \leftarrow \max \left\{ {0,{I_{t - 1}} + {P_t} - B_{t - 1}^{out} - D_t^{out}} \right\}$ \label{line-simu:It0}\;
$B_t^{out} \leftarrow  - \min \left\{ {0,{I_{t - 1}} + {P_t} - B_{t - 1}^{out} - D_t^{out}} \right\}$ \label{line-simu:Btout}\;
\lFor {each raw material item $i$}{
${P_{t + {l_{t,i}},i}} \leftarrow {O_{t,i}}$ \label{line-simu:Ptltii}
}
${r_t} \leftarrow \min \left\{ {\frac{{I_t^0}}{{\left( {{O_t} + O_{t - 1}^b} \right) \times {A^T}}},1} \right\}$ \label{line-simu:rt}\;
\For {each final product or sub-assembly item $i$}{
${M_{t,i}} \leftarrow \mathop {\min }\limits_{j,{a_{j,i}} \ne 0} \left\{ {{r_{t,j}}} \right\}\left( {{O_{t,i}} + O_{t - 1,i}^b} \right)$ \label{line-simu:Mti}\;
${P_{t + {l_{t,i}},i}} \leftarrow {M_{t,i}}$ \label{line-simu:Ptltii2}\;
}
$O_t^b \leftarrow {O_t} + O_{t - 1}^b - {M_t}$ \label{line-simu:Otb}\;
${I_t} \leftarrow I_t^0 - {M_t} \times {A^T}$ \label{line-simu:It}\;
${C_t} \leftarrow {I_t} \times {h^T} + B_t^{out} \times {p^T}$ \label{line-simu:Ct}\;
$C_t^{sum} \leftarrow C_{t - 1}^{sum} + {C_t}$  \label{line-simu:Ctsum}\;

}
\end{algorithm}

\subsection{Complexity Analysis}\label{subsec:complex-anal}

Computational complexity is often used to understand the efficiency of an algorithm and to compare the efficiency of different algorithms. In this paper we also adopt this tool to compare different algorithms. However, we want to emphasize that computational complexity does not completely explain the differences of the run times of different algorithms. For instance, based on our experience, tensorization may reduce the run time of our simulation algorithm by orders of magnitude, but it does not reduce the algorithm's complexity. Nevertheless, it provides a good measure to compare algorithms when all of them are properly tensorized.

The computational complexity of the addition/multiplication of two scalars is defined as ${\rm O}\left( 1 \right)$. Then, the computational complexity of vector (with $n$ entries) addition is ${\rm O}\left( n \right)$, and that of the multiplication of a $1 \times n$ vector and a $n \times n$ matrix is ${\rm O}\left( {{n^2}} \right)$. Therefore, in Algorithm \ref{alg:simulation}, the computational complexities of lines \ref{line-simu:Ptltii}, \ref{line-simu:Mti}, \ref{line-simu:Ptltii2}, \ref{line-simu:Ctsum} are ${\rm O}\left( 1 \right)$, those of lines \ref{line-simu:Ot1}, \ref{line-simu:Ot2}, \ref{line-simu:It0}, \ref{line-simu:Btout}, \ref{line-simu:Otb}, \ref{line-simu:Ct} are ${\rm O}\left( n \right)$, and those of lines \ref{line-simu:IPt}, \ref{line-simu:IPttemp}, \ref{line-simu:rt}, \ref{line-simu:It} are ${\rm O}\left( {{n^2}} \right)$, respectively. Furthermore, in each period, lines \ref{line-simu:Ptltii}, \ref{line-simu:Mti} and \ref{line-simu:Ptltii2} are executed ${\rm O}\left( n \right)$ times and all other parts are executed once. Given that there are totally $T$ periods in the simulation, the overall computational complexity of a simulation replication is ${\rm O}\left( {T{n^2}} \right)$.

If the matrix is sparse with density $\rho $, the computational complexity of  multiplication of a $1 \times n$ vector and a $n \times n$ sparse matrix reduces to ${\rm O}\left( {\rho {n^2}} \right)$. In our problem, the adjacency matrix is typically sparse. Then, utilizing sparse matrix techniques in the simulation algorithm not only reduces the memory consumption, but also improves the computing efficiency. So, the overall computational complexity reduces to ${\rm O}\left( {T\rho {n^2}} \right)$. If $n\rho\sim {\rm constant}$ as we showed in Equation (\ref{eqn:sparsity}), the overall computational complexity is ${\rm O}\left( {Tn} \right)$. We summarize these results in the following theorem.
\begin{theorem}
\label{theorem1}
For a general inventory system with $n$ items, the overall computational complexity of a single replication of Algorithm \ref{alg:simulation} with $T$ periods is ${\rm O}\left( {T{n^2}} \right)$. Using the sparse matrix techniques to handle the adjacency matrix, the complexity is ${\rm O}\left( {T\rho {n^2}} \right)$. Furthermore, if Equation (\ref{eqn:sparsity}) holds, the complexity is ${\rm O}\left( {Tn} \right)$.
\end{theorem}

\begin{remark}
Even though Theorem \ref{theorem1} is easy to derive, it has important implications. First, it showed that the use of sparse matrix techniques in the simulation algorithm can reduce the computational complexity from ${\rm O}\left( {T{n^2}} \right)$ to ${\rm O}\left( {Tn} \right)$. This is a significant reduction especially when handling large-scale inventory systems, where $n$ is in the orders of $10^4$ to $10^5$. Second, it provides a benchmark result to understand the performance of the gradient computation algorithms that are discussed in next section.
\end{remark}

\section{Simulation Optimization}\label{sec:opt}

Once we have the simulation model, our next step is to develop a gradient-based simulation optimization algorithm to solve the inventory optimization problem (\ref{eq:optprob}). To develop such an algorithm, we need to find an efficient algorithm to estimate the gradient of the expected total inventory cost {\small$E\left[{\sum\limits_{t = 1}^T {{C_t}\left( S \right)} }\right]$} with respect to the base-stock levels $S$.

\subsection{Gradient Computation}\label{subsec:grad-comput}

The simplest method for gradient estimation may be finite difference approximations. It is easy to understand and implement. However, it produces biased estimators and the computation time is prohibitively long, because at least $n+1$ simulation runs are necessary to compute just one observation of the gradient. Single-run unbiased gradient estimation methods have been proposed in the simulation literature, including the IPA method and the likelihood ratio method. Between these two methods, it is reported that when both are applicable, the IPA method typically produces gradient estimators that have smaller variances than the likelihood ratio method does \citep{glasserman2003monte}. 

\subsubsection{Infinitestimal Perturbation Analysis.}\label{subsubsec:ipa}

IPA is one of the most important gradient estimation methods in the field of stochastic simulation \citep{ho1983infinitesimal}. Based on \cite{glasserman2003monte}, we can estimate the derivative of $E\left[ {J\left( \theta  \right)} \right]$ using $E\left[ {\frac{d}{{d\theta }}J\left( \theta  \right)} \right]$, where ${\frac{d}{{d\theta }}J\left( \theta  \right)}$ is known as the IPA estimator (or pathwise gradient estimator), if  following conditions are satisfied:
\begin{itemize}
    \item ${J\left( \theta  \right)}$ is differentiable with probability 1, and
    \item ${J\left( \theta  \right)}$ is stochastically Lipschitz.
\end{itemize}

In our simulation model, the mathematical operations are typically addition, multiplication,  $\min \left( {} \right)$ or $\max \left( {} \right)$, which satisfy the above conditions. Besides, a detailed validation for using the IPA method in gradient estimation for multi-echelon production-inventory systems can be seen in \cite{glasserman1995sensitivity}.

Let $I_t^{temp} = {I_{t - 1}} + {P_t} - B_{t - 1}^{out} - D_t^{out}$, $I_t^{need} = \left( {{O_t} + O_{t - 1}^b} \right) \times {A^T}$ and ${k_{t,i}} = \mathop {\min }\limits_{j,{a_{j,i}} \ne 0} \left\{ {{r_{t,j}}} \right\}$, and let $E$ denote the identity matrix. Corresponding to each step in the simulation model,  the gradient of the total sample cost with respect to the base-stock levels can be calculated using the procedure shown in Algorithm \ref{alg:IPA}.
\begin{algorithm}
\small
\label{alg:IPA}
\SetAlgoNoLine
\caption{Procedure for Gradient Computation Using IPA}
\KwSty{Initialization:} $\frac{{\partial I{P_0}}}{{\partial S}} = \frac{{\partial {I_0}}}{{\partial S}} = 0,\frac{{\partial {O_0}}}{{\partial S}} = 0,\frac{{\partial B_0^{out}}}{{\partial S}} = 0,\frac{{\partial O_0^b}}{{\partial S}} = 0,\frac{{\partial C_0^{sum}}}{{\partial S}} = 0$\;
\For {$t = 1$ \KwTo $T$}{
$\frac{{\partial I{P_t}}}{{\partial S}} \leftarrow \frac{{\partial I{P_{t - 1}}}}{{\partial S}} + \frac{{\partial {O_{t - 1}}}}{{\partial S}} - A \times \frac{{\partial {O_{t - 1}}}}{{\partial S}}$ \label{line-ipa:IPt}\;
$\frac{{\partial {O_t}}}{{\partial S}} \leftarrow diag\left( { - {1_{\left\{ {\left( {I{P_t} - D_t^{out} - S} \right) < 0} \right\}}}} \right) \times \left( {\frac{{\partial I{P_t}}}{{\partial S}} - E} \right)$ \label{line-ipa:Ot1}\;
\For {$k = 1$ \KwTo ${n_l-1}$}{
$\frac{{\partial IP_t^{temp}}}{{\partial S}} \leftarrow \frac{{\partial I{P_t}}}{{\partial S}} - A \times \frac{{\partial {O_t}}}{{\partial S}}$ \label{line-ipa:IPttemp}\;
$\frac{{\partial {O_t}}}{{\partial S}} \leftarrow diag\left( { - {1_{\left\{ {\left( {IP_t^{temp} - S} \right) < 0} \right\}}}} \right) \times \left( {\frac{{\partial IP_t^{temp}}}{{\partial S}} - E} \right)$ \label{line-ipa:Ot2}\;
}

$\frac{{\partial I_t^{temp}}}{{\partial S}} \leftarrow \frac{{\partial {I_{t - 1}}}}{{\partial S}} + \frac{{\partial {P_t}}}{{\partial S}} - \frac{{\partial B_{t - 1}^{out}}}{{\partial S}}$ \label{line-ipa:Ittemp}\;
$\frac{{\partial I_t^0}}{{\partial S}} \leftarrow diag\left( {{1_{\left\{ {I_t^{temp} > 0} \right\}}}} \right) \times \frac{{\partial I_t^{temp}}}{{\partial S}}$\;
$\frac{{\partial B_t^{out}}}{{\partial S}} \leftarrow diag\left( { - {1_{\left\{ {I_t^{temp} \le 0} \right\}}}} \right) \times \frac{{\partial I_t^{temp}}}{{\partial S}}$ \label{line-ipa:Bt}\;

\lFor {each raw material item $i$}{
$\frac{{\partial {P_{t + {l_{t,i}},i}}}}{{\partial S}} \leftarrow \frac{{\partial {O_{t,i}}}}{{\partial S}}$ \label{line-ipa:Ptlii}
}

$\frac{{\partial I_t^{need}}}{{\partial S}} \leftarrow A \times \left( {\frac{{\partial {O_t}}}{{\partial S}} + \frac{{\partial O_{t - 1}^b}}{{\partial S}}} \right)$ \label{line-ipa:Itneed}\;

$\frac{{\partial {r_t}}}{{\partial S}} \leftarrow diag\left( {{1_{\left\{ {\frac{{I_t^0}}{{I_t^{need}}} < 1} \right\}}}} \right) \times diag\left( {\frac{1}{{I_t^{need}}}} \right) \times \left( {\frac{{\partial I_t^0}}{{\partial S}} - diag\left( {\frac{{I_t^0}}{{I_t^{need}}}} \right) \times \frac{{\partial I_t^{need}}}{{\partial S}}} \right)$ \label{line-ipa:rt}\;

\For {each final product or sub-assembly item $i$}{
$\frac{{\partial {k_{t,i}}}}{{\partial S}} \leftarrow {\left( {{1_{\left\{ {{r_{t,j}} = {k_{t,i}},{a_{j,i}} \ne 0} \right\}}}} \right)_{1 \times n}} \times \frac{{\partial {r_t}}}{{\partial S}}$ \label{line-ipa:kti}\;
$\frac{{\partial {M_{t,i}}}}{{\partial S}} \leftarrow \left( {{O_{t,i}} + O_{t - 1,i}^b} \right)\frac{{\partial {k_{t,i}}}}{{\partial S}} + {k_{t,i}}\left( {\frac{{\partial {O_{t,i}}}}{{\partial S}} + \frac{{\partial O_{t - 1,i}^b}}{{\partial S}}} \right)$\;
$\frac{{\partial {P_{t + {l_{t,i}},i}}}}{{\partial S}} \leftarrow \frac{{\partial {M_{t,i}}}}{{\partial S}}$ \label{line-ipa:Ftlii}\;
}

$\frac{{\partial {I_t}}}{{\partial S}} \leftarrow \frac{{\partial I_t^0}}{{\partial S}} - A \times \frac{{\partial {M_t}}}{{\partial S}}$ \label{line-ipa:It}\;

$\frac{{\partial O_t^b}}{{\partial S}} \leftarrow \frac{{\partial {O_t}}}{{\partial S}} + \frac{{\partial O_{t - 1}^b}}{{\partial S}} - \frac{{\partial {M_t}}}{{\partial S}}$ \label{line-ipa:Mtb}\;
$\frac{{\partial {C_t}}}{{\partial S}} \leftarrow h \times \frac{{\partial {I_t}}}{{\partial S}} + p \times \frac{{\partial B_t^{out}}}{{\partial S}}$ \label{line-ipa:Ct}\;

$\frac{{\partial C_t^{sum}}}{{\partial S}} \leftarrow \frac{{\partial C_{t - 1}^{sum}}}{{\partial S}} + \frac{{\partial {C_t}}}{{\partial S}}$ \label{line-ipa:Ctsum}\;
}
\end{algorithm}

It can be seen that most operations in Algorithm \ref{alg:IPA} are the operations of the Jacobian matrices. Given that the computational complexities of the addition and multiplication of two scalars are both ${\rm O}\left( 1 \right)$, the computational complexities of the addition and multiplication of two $n\times n$ matrices are ${\rm O}\left( {{n^2}} \right)$ and ${\rm O}\left( {{n^3}} \right)$, respectively. Therefore, for the steps in the procedure, the computation complexities of lines \ref{line-ipa:Ptlii}, \ref{line-ipa:kti}-\ref{line-ipa:Ftlii} are \ref{line-ipa:Ctsum} are all ${\rm O}\left( n \right)$, those of lines \ref{line-ipa:Ot1}, \ref{line-ipa:Ot2}-\ref{line-ipa:Bt}, \ref{line-ipa:rt}, \ref{line-ipa:Mtb} and \ref{line-ipa:Ct} are all ${\rm O}\left( {{n^2}} \right)$, and those of lines \ref{line-ipa:IPt}, \ref{line-ipa:IPttemp}, \ref{line-ipa:Itneed} and \ref{line-ipa:It} are all ${\rm O}\left( {{n^3}} \right)$. Furthermore, lines  \ref{line-ipa:Ptlii} and \ref{line-ipa:kti}-\ref{line-ipa:Ftlii} are all executed ${\rm O}\left( n \right)$ times and other lines are all executed ${\rm O}\left( 1 \right)$ times in each period. Hence, the overall computational complexity of calculating the IPA estimator over $T$ periods is ${\rm O}\left( {T{n^3}} \right)$. Considering the sparsity of the BOM network and Equation (\ref{eqn:sparsity}), we have the following theorem on the computational complexity of Algorithm \ref{alg:IPA}.
\begin{theorem}
\label{theorem2}
For a general inventory system with $n$ items, the computational complexity of Algorithm \ref{alg:IPA} for a simulation with $T$ periods is ${\rm O}\left( {T{n^3}} \right)$. Suppose that the sparse matrix techniques are used to handle the adjacency matrix, the computational complexity is ${\rm O}\left( {T\rho {n^3}} \right)$. Furthermore, if Equation (\ref{eqn:sparsity}) holds, the computational complexity may be reduced to ${\rm O}\left( {T{n^2}} \right)$.
\end{theorem}

\begin{remark}
\label{remark.finite}
Recall that the computational complexity of the simulation algorithm is approximately ${\rm O}\left( {Tn} \right)$ when Equation (\ref{eqn:sparsity}) holds, the computational complexity of calculating the IPA estimator is actually of the same order as the finite difference estimation, which requires running the simulation model at least $n + 1$ times. This result seems somehow counter-intuitive, because we typically think the IPA estimator, which requires only a single replication, is computationally more efficient than the finite-difference estimators, which require $O(n)$ replications. The reason is that, although the IPA calculates the sample path gradient in just one simulation, its calculation requires repeatedly handling of the Jacobian matrices, which requires far more computing time than the vector computations in the simulation algorithm. 
\end{remark}

Numerical results also show that for large-scale inventory systems, gradient estimation using IPA is rather time-consuming. In order to meet the needs of further simulation optimization, more efficient gradient estimation methods need to be considered.

\subsubsection{Back Propagation.}\label{subsubsec:bp}

As we have shown in Section \ref{sec:rnn}, the training of large-scale RNNs also need to compute the sample-path gradient of the sample total loss with respect to the weights of the neural networks, which is typically done by the BP algorithms. The BP algorithm for training neural networks was proposed by \cite{rumelhart1986learning}, and it is one of the corner stones of deep learning. The idea of the BP algorithm is to compute the gradient using a reverse mode, after finishing the calculation of the function value. It is different from the IPA algorithm, which uses a forward mode to calculate the gradient alongside the calculation of the function value. It has been shown that, when the input is of a multi-dimensional vector (such as the weight vector in RNN training or the base-stock level vector in inventory optimization), the BP algorithm is typically computationally more efficient than the forward mode algorithms (see, for instance, \cite{giles2006smoking}).

\begin{figure}
\begin{center}
\includegraphics[height=2in]{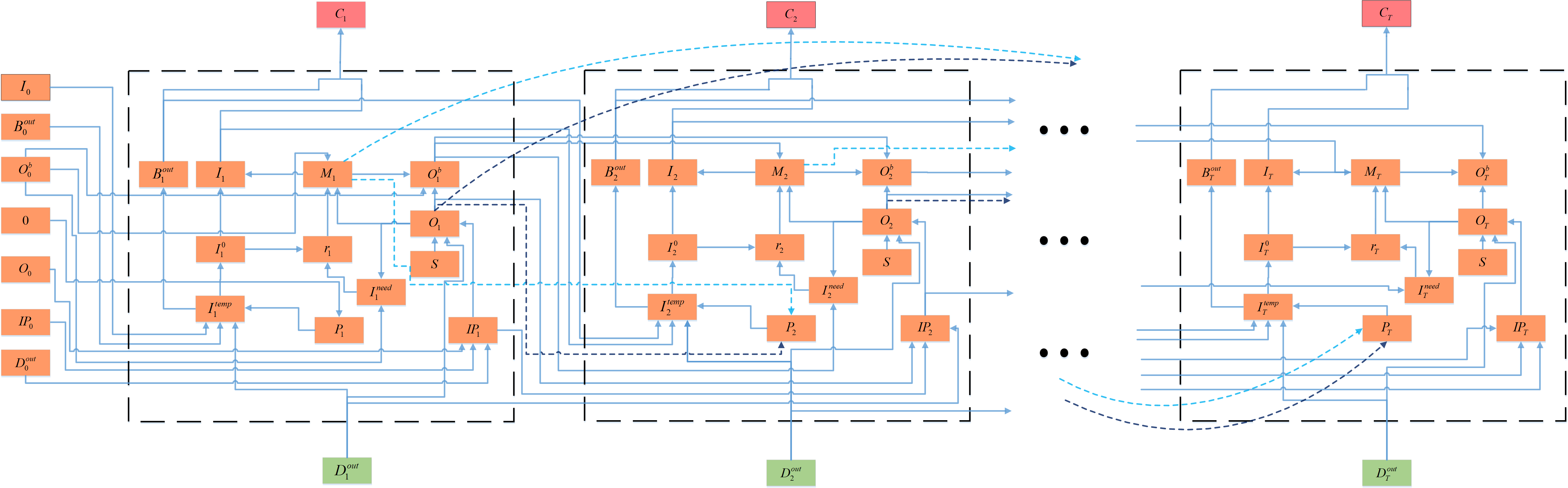}
\caption{Computational Graph}
\label{fig:comp graph}
\end{center}
\end{figure}

In deep learning frameworks like TensorFlow, PyTorch, etc., BP algorithms are typically implemented through computational graphs. A computational graph is a directed graph for expressing and evaluating mathematical expressions of an algorithm. We construct the computational graph of Algorithm \ref{alg:simulation} and present it in Figure \ref{fig:comp graph}. It shows how the simulation algorithm is executed, and it serves as a map for us to develop the BP algorithm for the inventory model.

Notice that the BP algorithm calculates the gradient after finishing the calculation of the function value. In our inventory simulation model, it requires recording certain gradient-related information in the process of simulation before entering in the backward stage. We summarize this part in Algorithm \ref{alg:bp-part1}.
\begin{algorithm}
\small
\label{alg:bp-part1}
\SetAlgoNoLine
\caption{Procedure for Gradient Computation Using BP Part $1$}
\For {$t = 1$ \KwTo $T$}{
push $\frac{{\partial {O_t}}}{{\partial I{P_t}}} = diag\left( { - {1_{\left\{ {\left( {I{P_t} - D_t^{out} - S} \right) < 0} \right\}}}} \right)$ into Stack$1$\;
\lFor{$k = 1$ \KwTo ${n_l}-1$}{
push $\frac{{\partial {O_t}}}{{\partial IP_t^{temp}}} = diag\left( { - {1_{\left\{ {\left( {IP_t^{temp} - S} \right) < 0} \right\}}}} \right)$ into Stack$1$}
record $\frac{{\partial I_t^0}}{{\partial I_t^{temp}}} = diag\left( {{1_{\left\{ {I_t^{temp} > 0} \right\}}}} \right)$ and $\frac{{\partial B_t^{out}}}{{\partial I_t^{temp}}} = diag\left( { - {1_{\left\{ {I_t^{temp} \le 0} \right\}}}} \right)$\;

\For {each item $i$ which has procurement order in period $t$}{
record $i$ in List$1$ of period ${t + {l_{t,i}}}$\;
record ${l_{t,i}}$ in List$2$ of period ${t + {l_{t,i}}}$\;
}

record $\frac{{\partial {r_t}}}{{\partial I_t^0}} = diag\left( {{1_{\left\{ {\frac{{I_t^0}}{{I_t^{need}}} < 1} \right\}}}} \right) \times diag\left( {\frac{1}{{I_t^{need}}}} \right)$ and $\frac{{\partial {r_t}}}{{\partial I_t^{need}}} = diag\left( {{1_{\left\{ {\frac{{I_t^0}}{{I_t^{need}}} < 1} \right\}}}} \right) \times diag\left( {\frac{1}{{I_t^{need}}}} \right) \times diag\left( {\frac{{I_t^0}}{{I_t^{need}}}} \right)$\;

\For {each item $i$ whose started production is nonzero in period $t$}{
record $i$ in List$3$ of period ${t + {l_{t,i}}}$\;
record ${l_{t,i}}$ in List$4$ of period ${t + {l_{t,i}}}$\;
record ${{k_{t,i}}}$ and $\frac{{\partial {M_{t,i}}}}{{\partial {r_t}}} = \left( {{O_{t,i}} + O_{t - 1,i}^b} \right){\left( {{1_{\left\{ {{r_{t,j}} = {k_{t,i}},{a_{j,i}} \ne 0} \right\}}}} \right)_{1 \times n}}$\;
}
}
\end{algorithm}

After recording gradient-related information during the simulation process, the gradient of the sample total cost with respect to the base-stock levels may be computed form the last period $T$ backward to period $1$. The procedure is summarized in Algorithm \ref{alg:bp-part2}.
\begin{algorithm}
\small
\label{alg:bp-part2}
\SetAlgoNoLine
\caption{Procedure for Gradient Computation Using BP Part $2$}
\KwSty{Initialization:} $\frac{{\partial C_T^{sum}}}{{\partial I_{T + 1}^{temp}}} = 0,\frac{{\partial C_T^{sum}}}{{\partial I{P_T}}} = 0,\frac{{\partial C_T^{sum}}}{{\partial S}} = 0,\frac{{\partial C_T^{sum}}}{{\partial O_T^b}} = 0,\frac{{\partial C_T^{sum}}}{{\partial {r_t}}} = 0,\frac{{\partial C_T^{sum}}}{{\partial {O_t}}} = 0,\frac{{\partial C_T^{sum}}}{{\partial {M_t}}} = 0\left( {t = 1, \ldots ,T} \right)$\;
\For {$t = T$ \KwTo $1$}{
$\frac{{\partial C_T^{sum}}}{{\partial {I_t}}} \leftarrow h + \frac{{\partial C_T^{sum}}}{{\partial I_{t + 1}^{temp}}}$,
$\frac{{\partial C_T^{sum}}}{{\partial B_t^{out}}} \leftarrow p - \frac{{\partial C_T^{sum}}}{{\partial I_{t + 1}^{temp}}}$ \label{line-bp2:It1}\;

$\frac{{\partial C_T^{sum}}}{{\partial {M_t}}} \leftarrow \frac{{\partial C_T^{sum}}}{{\partial {M_t}}} - \frac{{\partial C_T^{sum}}}{{\partial {I_t}}} \times A - \frac{{\partial C_T^{sum}}}{{\partial O_t^b}}$ \label{line-bp2:Mt}\;

$\frac{{\partial C_T^{sum}}}{{\partial {{O'}_t}}} \leftarrow \frac{{\partial C_T^{sum}}}{{\partial O_t^b}} + \frac{{\partial C_T^{sum}}}{{\partial {M_t}}} \times diag\left( {{{\left( {{k_{t,i}}} \right)}_{1 \times n}}} \right)$ \label{line-bp2:Ot'}\;

\For {each item $i$ whose started production is nonzero in period $t$}{
$\frac{{\partial C_T^{sum}}}{{\partial {r_t}}} \leftarrow \frac{{\partial C_T^{sum}}}{{\partial {r_t}}} + \frac{{\partial C_T^{sum}}}{{\partial {M_{t,i}}}}\frac{{\partial {M_{t,i}}}}{{\partial {r_t}}}$ \label{line-bp2:rt}\;
}

$\frac{{\partial C_T^{sum}}}{{\partial I_t^0}} \leftarrow \frac{{\partial C_T^{sum}}}{{\partial {I_t}}} + \frac{{\partial C_T^{sum}}}{{\partial {r_t}}} \times \frac{{\partial {r_t}}}{{\partial I_t^0}}$ \label{line-bp2:It0}\; 

$\frac{{\partial C_T^{sum}}}{{\partial I_t^{need}}} \leftarrow \frac{{\partial C_T^{sum}}}{{\partial {r_t}}} \times \frac{{\partial {r_t}}}{{\partial I_t^{need}}}$ \label{line-bp2:Itneed}\;
$\frac{{\partial C_T^{sum}}}{{\partial {{O'}_t}}} \leftarrow \frac{{\partial C_T^{sum}}}{{\partial {{O'}_t}}} + \frac{{\partial C_T^{sum}}}{{\partial I_t^{need}}} \times A$\;
$\frac{{\partial C_T^{sum}}}{{\partial I_t^{temp}}} \leftarrow \frac{{\partial C_T^{sum}}}{{\partial I_t^0}} \times \frac{{\partial I_t^0}}{{\partial I_t^{temp}}} + \frac{{\partial C_T^{sum}}}{{\partial B_t^{out}}} \times \frac{{\partial B_t^{out}}}{{\partial I_t^{temp}}}$\; 

$\frac{{\partial C_T^{sum}}}{{\partial {O_t}}} \leftarrow \frac{{\partial C_T^{sum}}}{{\partial {O_t}}} + \frac{{\partial C_T^{sum}}}{{\partial {{O'}_t}}} + \frac{{\partial C_T^{sum}}}{{\partial I{P_t}}} \times \left( {E - A} \right)$ \label{line-bp2:Ot}\;
\For{$k = 1$ \KwTo ${n_l}-1$}{
pop $\frac{{\partial {O_t}}}{{\partial IP_t^{temp}}}$ from Stack$1$\;
$\frac{{\partial C_T^{sum}}}{{\partial IP_t^{temp}}} \leftarrow \frac{{\partial C_T^{sum}}}{{\partial {O_t}}} \times \frac{{\partial {O_t}}}{{\partial IP_t^{temp}}}$,
$\frac{{\partial C_T^{sum}}}{{\partial S}} \leftarrow \frac{{\partial C_T^{sum}}}{{\partial S}} - \frac{{\partial C_T^{sum}}}{{\partial IP_t^{temp}}}$ \label{line-bp2:IPttemp}\;
$\frac{{\partial C_T^{sum}}}{{\partial I{P_t}}} \leftarrow \frac{{\partial C_T^{sum}}}{{\partial I{P_t}}} + \frac{{\partial C_T^{sum}}}{{\partial IP_t^{temp}}}$,
$\frac{{\partial C_T^{sum}}}{{\partial {O_t}}} \leftarrow  - \frac{{\partial C_T^{sum}}}{{\partial IP_t^{temp}}} \times A$  \label{line-bp2:IPt-1}\;
}

pop $\frac{{\partial {O_t}}}{{\partial I{P_{t}}}}$ from Stack$1$\; 

$\frac{{\partial C_T^{sum}}}{{\partial S}} \leftarrow \frac{{\partial C_T^{sum}}}{{\partial S}} - \frac{{\partial C_T^{sum}}}{{\partial {O_t}}} \times \frac{{\partial {O_t}}}{{\partial I{P_t}}}$ \label{line-bp2:S}\;

$\frac{{\partial C_T^{sum}}}{{\partial I{P_{t - 1}}}} \leftarrow \frac{{\partial C_T^{sum}}}{{\partial I{P_t}}} + \frac{{\partial C_T^{sum}}}{{\partial {O_t}}} \times \frac{{\partial {O_t}}}{{\partial I{P_t}}}$ \label{line-bp2:IPt-1-2}\;

$\frac{{\partial C_T^{sum}}}{{\partial O_{t - 1}^b}} \leftarrow \frac{{\partial C_T^{sum}}}{{\partial {{O'}_t}}}$,
$\frac{{\partial C_T^{sum}}}{{\partial {P_t}}} \leftarrow \frac{{\partial C_T^{sum}}}{{\partial I_t^{temp}}}$ \label{line-bp2:Ot-1b}\;

\For {each item $i$ in List$1$}{
get the recorded lead time ${l_{u,i}}$ in List$2$ ($u = t - {l_{u,i}}$) \label{line-bp2:list2}\;
$\frac{{\partial C_T^{sum}}}{{\partial {O_{t - {l_{u,i}},i}}}} \leftarrow \frac{{\partial C_T^{sum}}}{{\partial {P_{t,i}}}}$ \label{line-bp2:Ot-luii}\;
}
\For {each item $i$ in List$3$}{
get the recorded lead time ${l_{u,i}}$ in List$4$ ($u = t - {l_{u,i}}$) \label{line-bp2:list4}\;
$\frac{{\partial C_T^{sum}}}{{\partial {M_{t - {l_{u,i}},i}}}} \leftarrow \frac{{\partial C_T^{sum}}}{{\partial {P_{t,i}}}}$ \label{line-bp2:Mt-luii}\;
}
}
\end{algorithm}

Different from the IPA algorithm (Algorithm \ref{alg:IPA}), which is mainly composed of operations of Jacobian matrices, the BP algorithm (Algorithms \ref{alg:bp-part1} and \ref{alg:bp-part2}) is mainly composed of vector operations. In Algorithm \ref{alg:bp-part1}, the computational complexities of lines 2-4 and 8 are ${\rm O}\left( n \right)$, and those of lines 6-7 and 10-12 are ${\rm O}\left( 1 \right)$ and they are executed ${\rm O}\left( n \right)$ times in each period. Therefore, the computational complexity of Algorithm \ref{alg:bp-part1} is ${\rm O}\left( {Tn} \right)$.  In Algorithm \ref{alg:bp-part2}, the computational complexities of lines \ref{line-bp2:list2}-\ref{line-bp2:Ot-luii} and \ref{line-bp2:list4}-\ref{line-bp2:Mt-luii} are all ${\rm O}\left( 1 \right)$, those of lines \ref{line-bp2:It1}, \ref{line-bp2:Ot'}, \ref{line-bp2:rt} and \ref{line-bp2:Ot-1b} are all ${\rm O}\left( n \right)$, and those of lines \ref{line-bp2:Mt}, \ref{line-bp2:It0}-\ref{line-bp2:Ot}, \ref{line-bp2:IPttemp}-\ref{line-bp2:IPt-1} and \ref{line-bp2:S}-\ref{line-bp2:IPt-1-2} are all ${\rm O}\left( {{n^2}} \right)$. Because lines \ref{line-bp2:rt}, \ref{line-bp2:list2}-\ref{line-bp2:Ot-luii} and \ref{line-bp2:list4}-\ref{line-bp2:Mt-luii} are executed ${\rm O}\left( n \right)$ times, while others are executed ${\rm O}\left( 1 \right)$ times in each time period, the computational complexity of Algorithm \ref{alg:bp-part2} is therefore ${\rm O}\left( {T{n^2}} \right)$. Then, we have the following theorem on the computational complexity of the BP algorithm.
\begin{theorem}
\label{theorem3}
For a general inventory system with $n$ items, the computational complexity of the BP algorithm (i.e., Algorithms \ref{alg:bp-part1} and \ref{alg:bp-part2}) for a simulation with $T$ periods is ${\rm O}\left( {T{n^2}} \right)$. Suppose that the sparse matrix techniques are used to handle the adjacency matrix, the computational complexity is ${\rm O}\left( {T\rho {n^2}} \right)$. Furthermore, if Equation (\ref{eqn:sparsity}) holds, the computational complexity may be reduced to ${\rm O}\left( {Tn} \right)$.
\end{theorem}

With the BP algorithm, the computational complexity of gradient calculation is now the same as that of the simulation algorithm and it is of an order $n$ faster than the IPA algorithm. In the numerical results reported in Section \ref{sec:exp}, we find that the BP algorithm and IPA algorithm produce the same gradient estimates (as expected), but the BP algorithm is significantly faster especially when $n$ is large.

\subsection{Optimization Algorithm}\label{subsec:opt-algo}

As described in Section \ref{sec:state}, the objective of the inventory optimization problem is to find the optimal base-stock levels $S$ that minimize the expected cumulative cost over $T$ periods, i.e.,
\[\min\quad E\left[ {\sum\limits_{t = 1}^T {{C_t}\left( S \right)} } \right].\]
To solve this, we can compute the sample-path gradient with respect to the base-stock levels and apply SGD algorithms. However, by applying SGD algorithms directly, we observe that the solution typically keeps inventory at almost all nodes and it is very difficult to implement in practice. Instead, we want to find solutions that keep inventory only at a small fraction of the nodes. To solve the problem, we add the $L_1$ norm of the decision variables along with a tuning parameter $\lambda>0$, and solve the following problem:
\begin{equation}\label{eq:L1}
    \min\quad {E\left[ {\sum\limits_{t = 1}^T {{C_t}\left( S \right)} } \right] + \lambda {{\left\| S \right\|}_1}} .
\end{equation}
The $L_1$ regularization is also known as Lasso in the statistics literature \citep{tibshirani1996regression}, and it is known to introduce  sparsity to the solution.

From another perspective, for large-scale inventory optimization problems, the number of samples used in the optimization process is not very large, typically less than number of decision variables (i.e., the base-stock levels). Thus, it is an over-parameterized model \citep{soltanolkotabi2018theoretical}, and the bias of the sample optimal value cannot be ignored. To reduce the bias in the optimal solution, regularization methods are also typically used \citep{buhlmann2011statistics}. Therefore, the $L_1$ regularization introduced in Problem (\ref{eq:L1}) not only keeps the inventory only at a small fraction of the nodes, it also helps reduce the bias.

Notice that the $L_1$ norm term in Problem (\ref{eq:L1}) is not smooth. One way to solve the optimization problem is to use the stochastic subgradient descent (SSGD) algorithm \citep{shor2012minimization}. Let {\small$F\left( S \right) = E\left[ {f\left( S \right)} \right] = E\left[ {\sum\limits_{t = 1}^T {{C_t}\left( S \right)} } \right]$} and $G\left( S \right) = \lambda {\left\| S \right\|_1}$. The iteration of the SSGD algorithm takes the form
\[{S_{k + 1}} = {S_k} - {t_k}\left( {\nabla f\left( {{S_k}} \right) + {\xi _k}} \right),\]
where ${t_k}$ is the step-size and ${{\xi _k}}$ denotes the subgradient of $G$ at $S_k$. A main drawback of this method is its lack of capability in exploiting the problem structure of the $L_1$ norm, thus having poor convergence properties \citep{xiao2009dual}. 

In contrast, there exist other optimization methods that exploit the problem structures and are better suited for this problem. Notice that the objective function of Problem (\ref{eq:L1}) is the sum of a smooth function and a simple convex non-smooth function. The proximal gradient method, sometimes called ISTA (iterative shrinkage-thresholding algorithm), can solve this type of problems, where the objective is a sum of a differentiable term and a non-differentiable convex function, with a convergence rate of ${\rm O}\left( {1/k} \right)$ \citep{parikh2014proximal}. By utilizing Nesterov acceleration method \citep{nesterov1983method}, FISTA (fast ISTA) promotes the convergence rate to ${\rm O}\left( {1/{k^2}} \right)$ \citep{beck2009fast}. FISTA is currently a popular algorithm in the field of machine learning. 

For optimization problem with objective
\[\min\ \ F\left( x \right) + G\left( x \right)\]
where $F$ is a smooth cost function, and $G$ is a possibly non-smooth regularization term, the basic iteration of FISTA is
\[{x_{k + 1}} = pro{x_{{t_{k + 1}}G}}\left( {{y_k} - {t_{k + 1}}\nabla F({y_k})} \right)\]
\[{y_{k + 1}} = {x_{k + 1}} + \frac{k}{{k + r}}\left( {{x_{k + 1}} - {x_k}} \right)\]
where $r \ge 3$, $\left\{ {{t_k}} \right\}$ is a sequence of positive and non-increasing step sizes, and $pro{x_{{t_{k + 1}}G}}\left(  \cdot  \right)$ is the proximal operator. If $G = \lambda {\left\| x \right\|_1}$, 
\[pro{x_{{t_{k + 1}}G}}{\left( x \right)_i} = sign\left( {{x_i}} \right)\max \left( {0,\left( {\left| {{x_i}} \right| - {t_{k + 1}}\lambda } \right)} \right).\]
When the function $F$ is represented as an expectation $F\left( x \right) = {E_\xi }\left[ {f\left( {\xi ,x} \right)} \right]$, as in our problem, we need to use a stochastic version of FISTA \citep{atchade2017perturbed,salimperformance}. The true gradient ${\nabla F({y_k})}$ is replaced by a sample-path gradient $\nabla f({y_k})$ in each iteration. When used with a constant step size, the stochastic FISTA achieves a convergence rate of {\small${\rm O}\left( {1/\sqrt k } \right)$}, and with decreasing step sizes it achieves a rate of {\small${\rm O}\left( {\log \left( k \right)/\sqrt k } \right)$}. Moreover, the stochastic FISTA is close to its ${\rm O}\left( {1/{k^2}} \right)$ deterministic behavior in the first several iterations \citep{salimperformance}. 

It is worth noting that the inventory optimization problem is different from the training of RNNs or other machine learning problems. Affected by the regularization term, the optimization solution may not be optimal for the original problem, i.e., Problem (\ref{eq:optprob}). Thus, a re-optimization step may be added to further improve the performance of the inventory decisions. Then, we have the following two-stage procedure.

\vspace{11pt}
\begin{descr}
\item[\textbf{Stage 1.}] 
Solve the optimization problem with the ${L_1}$ norm regularizer to select the inventory locations.
{\small
\[S_1^* = \arg \min \left\{ {E\left[ {\sum\limits_{t = 1}^T {{C_t}\left( S \right)} } \right] + \lambda {{\left\| S \right\|}_1}} \right\}.\]}
\item[\textbf{Stage 2.}] 
Based on the result from the previous step, fix certain
elements of $S$ to be $0$ and solve the modified problem:
{\small
\[\min\ E\left[ {\sum\limits_{t = 1}^T {{C_t}\left( \tilde S \right)} } \right]\quad{\rm where}\ \tilde S = S \odot {1_{\left\{ {S_1^* > 0} \right\}}},\]}
where $\odot$ denotes element-wise product.
\end{descr}

In Stage $1$, our goal is to select a small fraction of the nodes to keep inventory. In Stage $2$, we only consider the selected inventory locations and solve a much smaller-scale smooth inventory optimization problems to determine their base-stock levels. We use the stochastic FISTA in Stage 1 and the SGD in Stage 2.

Notice that solving the regularized problem in two stages is not new. In the field of high dimensional regression, \cite{meinshausen2007relaxed} proposed a two-stage procedure, termed the relaxed Lasso for high-dimensional data where the number of predictor variables $p$ is much larger than the number of observations $n$. Their theoretical and numerical results demonstrate that the two-stage procedure performs better than the regular Lasso estimator for high-dimensional data.

\section{Numerical Experiments}\label{sec:exp}

In this section, we first test the performances of the simulation and gradient computation algorithms on inventory systems of different scales. Then, we conduct a series of experiments to understand the behaviors of the optimization algorithm and to compare it with the GS model of \cite{graves2000optimizing} on small- to medium-scale problems where the GS model applies. In this section, all the computer programs are coded in Python and all the experiments are run on a computer with two Intel Xeon Gold 6248R CPUs (each with 24 cores) and 256GB RAM.

\subsection{Simulation and Gradient Computation }\label{subsec:performance-simu&grad}

To test the performance of our algorithms, We build test inventory examples based on the BOM characteristics that we observe from our consulting experience. We set the average degrees $\left\langle k \right\rangle $ of all BOM networks as $10$ and the time horizons as $100$ periods.

\subsubsection{Performance of the Simulation Algorithms.}

We test the traditional simulation algorithm, the tensorized simulation algorithms (i.e., Algorithm \ref{alg:simulation}) with dense matrices or sparse matrices for inventory systems with different number of nodes, ranging from medium-scale system with 1,000 nodes to very large-scale system with 500,000 nodes, and report the simulation run times of different algorithms, averaging over 10 independent replications, in Table \ref{tab:time consum-simu}, where ``-'' indicates that the algorithm takes over a day so it is terminated before completion and ``/'' indicates that there is not enough memory in our computer to run the algorithm.

\begin{table}[ht]\small
  \caption{Run Time of Simulation for a Single Replication}
  \label{tab:time consum-simu}
  \centering
  \renewcommand{\arraystretch}{1.1}
  \begin{tabular}{cccc}
  \toprule
    \multicolumn{1}{c}{\begin{tabular}[c]{@{}c@{}}Number\\  of nodes\end{tabular}} & \multicolumn{1}{c}{\begin{tabular}[c]{@{}c@{}}Traditional\\ algorithm\end{tabular}} & \multicolumn{1}{c}{\begin{tabular}[c]{@{}c@{}}Algorithm 1\\ dense \end{tabular}} & \multicolumn{1}{c}{\begin{tabular}[c]{@{}c@{}}Algorithm 1\\ sparse \end{tabular}}\\
   \midrule
    1000 & 2.21min & 0.28s & 0.19s\\
    5000 & 56.00min & 6.28s & 0.73s\\
    10000 & 212.18min & 25.12s & 1.48s\\
    50000 & - & 632.21s & 6.92s\\
    100000 & - & 2535.56s & 13.90s\\
    500000 & - & / & 85.70s\\
    \bottomrule
  \end{tabular}
\end{table}

There are several interesting findings from the results. First, the computational complexities of the traditional algorithm and the tensorized algorithm with dense matrices are approximately $O(n^2)$, while that of the tensorized algorithm with sparse matrices is approximately $O(n)$, which are consistent with Theorem \ref{theorem1}. Second, even though tensorization does not improve the computational complexity, it reduces the computational time significantly. Third, the use of sparse matrices reduces significantly both the memory requirements and the run times, further allowing the algorithm to handle very large-scale problems.

\subsubsection{Performance of the Gradient Computation Algorithms.}

We use the same test problems to test the two gradient computation algorithms proposed in this paper, the IPA algorithm (Algorithm \ref{alg:IPA}) and the BP Algorithm (Algorithms \ref{alg:bp-part1} and \ref{alg:bp-part2}). Notice that both algorithms are tensorized and both may use dense or sparse adjacency matrices. The run times of different algorithms are reported in Table \ref{tab:time consum-grad}.

\begin{table}[ht]\small
  \caption{Run Time of Gradient Computation for a Single Replication}
  \label{tab:time consum-grad}
  \centering
  \renewcommand{\arraystretch}{1.1}
  \begin{tabular}{ccccc}
  \toprule
   \multicolumn{1}{c}{\begin{tabular}[c]{@{}c@{}}Number\\  of nodes\end{tabular}} & IPA-dense & IPA-sparse & BP-dense & BP-sparse\\
   \midrule
    1000 & 0.55min & 0.34min & 1.00s & 0.81s \\
    5000 & 29.36min & 3.01min & 14.36s & 2.95s \\
    10000 &219.73min & 11.55min & 53.35s & 5.68s \\
    50000 & - & 654.67min & 1277.84s & 26.88s \\
    100000 & - & - & 5152.54s & 54.97s \\
    500000 & / & - & / & 305.64s \\
    \bottomrule
  \end{tabular}
\end{table}

From the results we see that the computational complexities of the IPA algorithms, with dense and sparse adjacency matrices, and the BP algorithms, with dense and sparse adjacency matrices, are approximately ${\rm O}\left( {{n^3}} \right)$, ${\rm O}\left( {{n^2}} \right)$, ${\rm O}\left( {{n^2}} \right)$ and ${\rm O}\left( {{n}} \right)$ respectively, which are consistent with Theorems \ref{theorem2} and \ref{theorem3}. Furthermore, it is clear that the BP algorithms significantly outperform the IPA algorithms in general, and the BP algorithm with sparse adjacency matrices is capable of handling very large-scale problems.

\subsubsection{Simulation and Gradient Computation Using TensorFlow.}

One major impediment to the use of our algorithms in general inventory optimization is the construction of the computational graph and the derivation of the sample-path gradient. Small changes of the logic of the simulation process may result in quite different computational graphs and thus very different gradient estimators. Therefore, the use of these algorithms may require a significant amount of analyst's effort, which may be very expensive in practice. This motivates us to program the tensorized simulation algorithm using the machine-learning tool TensorFlow, which may construct the computational graph and compute the sample-path gradient automatically. Moreover, the use of TensorFlow also makes it possible to utilize GPUs for matrix operations used in the algorithms, which may provide further speedups.

We implement the tensorized simulation algorithm using TensorFlow, and also use its automatic BP tool to compute the sample-path gradient. We test both algorithms with dense or sparse adjacency matrices. We further test the algorithms using an additional NVIDIA GeForce RTX 3090 GPU with 24GB memory. The results are reported in Tables \ref{tab:time consum-simu-tf} and \ref{tab:time consum-grad-tf}.

\begin{table}[ht]\small
  \caption{Run Time of Simulation Using TensorFlow for a Single Replication}
  \label{tab:time consum-simu-tf}
  \centering
  \renewcommand{\arraystretch}{1.1}
  \begin{tabular}{ccccc}
  \toprule
   \multicolumn{1}{c}{\begin{tabular}[c]{@{}c@{}}Number\\  of nodes\end{tabular}} & \multicolumn{1}{c}{\begin{tabular}[c]{@{}c@{}}TensorFlow-\\CPU-dense\end{tabular}} & \multicolumn{1}{c}{\begin{tabular}[c]{@{}c@{}}TensorFlow-\\GPU-dense\end{tabular}} & \multicolumn{1}{c}{\begin{tabular}[c]{@{}c@{}}TensorFlow-\\CPU-sparse\end{tabular}} & \multicolumn{1}{c}{\begin{tabular}[c]{@{}c@{}}TensorFlow-\\GPU-sparse\end{tabular}}\\
   \midrule
    1000 & 1.66s & 4.60s & 1.69s & 2.98s \\
    5000 & 13.08s & 12.76s & 7.08s & 11.61s \\
    10000 & 40.19s & 24.45s & 14.38s & 23.14s \\
    50000 & 735.64s & 129.72s & 71.30s & 112.70s \\
    100000 & 2912.50s & / & 145.31s & 231.39s \\
    500000 & / & / & 775.07s & 1144.78s \\ 
    \bottomrule
  \end{tabular}
\end{table}

\begin{table}[ht]\small
  \caption{Run Time of Gradient Computation Using TensorFlow for a Single Replication}
  \label{tab:time consum-grad-tf}
  \centering
  \renewcommand{\arraystretch}{1.1}
  \begin{tabular}{ccccc}
  \toprule
   \multicolumn{1}{c}{\begin{tabular}[c]{@{}c@{}}Number\\  of nodes\end{tabular}} & \multicolumn{1}{c}{\begin{tabular}[c]{@{}c@{}}TensorFlow-\\CPU-dense\end{tabular}} & \multicolumn{1}{c}{\begin{tabular}[c]{@{}c@{}}TensorFlow-\\GPU-dense\end{tabular}} & \multicolumn{1}{c}{\begin{tabular}[c]{@{}c@{}}TensorFlow-\\CPU-sparse\end{tabular}} & \multicolumn{1}{c}{\begin{tabular}[c]{@{}c@{}}TensorFlow-\\GPU-sparse\end{tabular}}\\
   \midrule
    1000 & 9.79s & 12.33s & 10.51s & 13.37s \\
    5000 & 59.94s & 54.59s & 48.06s & 55.02s \\
    10000 & 152.76s & 108.37s & 100.99s & 110.11s \\
    50000 & 2453.98s & / & 609.78s & 547.05s \\
    100000 & 11112.54s & / & 1874.00s & / \\
    500000 & / & / & 20337.55s & / \\ 
    \bottomrule
  \end{tabular}
\end{table}

There are several interesting findings from these results. First, when the inventory problems are of medium to large scales (e.g., 1,000 to 50,000 nodes), the TensorFlow implementations of the simulation algorithm and the automatic BP algorithm, both with sparse adjacency matrices, provide acceptable run times, making TensorFlow a very promising tool for inventory optimization. Second, it appears that the use of GPU does not provide any advantage when the sparse adjacency matrices are used, in terms of the computational time. Third, compared to Tables \ref{tab:time consum-simu} and \ref{tab:time consum-grad}, it is clear that our original method with sparse adjacency matrices is significantly faster than the TensorFlow implementations, making it an ideal tool to solve large- to very large-scale problems, though it requires a significant amount of analyst's effort to construct the computational graph and to derive the BP gradient.

\subsubsection{Performance of Multiple Replications.}

To solve the inventory optimization problems, the objective value and the gradient need to be evaluated based on multiple replications (i.e., a mini-batch) in each epoch (i.e., iteration of the optimization algorithms). Therefore, it is important to understand the run times of the algorithms when multiple replications are run in parallel. We test the algorithms with a mini-batch of 10 and report the results in Tables \ref{tab:time consum-simu-multi} and \ref{tab:time consum-grad-multi}.
\begin{table}[ht]\small
  \caption{Run Time of Simulation for Multiple Replications}
  \label{tab:time consum-simu-multi}
  \centering
  \renewcommand{\arraystretch}{1.1}
  \begin{tabular}{ccc}
  \toprule
   \multicolumn{1}{c}{\begin{tabular}[c]{@{}c@{}}Number\\  of nodes\end{tabular}} & \multicolumn{1}{c}{\begin{tabular}[c]{@{}c@{}}Single\\ replication\end{tabular}} & \multicolumn{1}{c}{\begin{tabular}[c]{@{}c@{}}Mini-batch\\ of 10\end{tabular}}\\
   \midrule
    1000 & 0.19s & 0.72s\\
    5000 & 0.73s & 1.75s\\
    10000 & 1.48s & 2.99s\\
    50000 & 6.91s & 13.47s\\
    100000 & 13.90s & 26.36s\\
    500000 & 86.42s & 149.63s\\ 
    \bottomrule
  \end{tabular}
\end{table}

\begin{table}[ht]\small
  \caption{Run Time of BP for Multiple Replications}
  \label{tab:time consum-grad-multi}
  \centering
  \renewcommand{\arraystretch}{1.1}
  \begin{tabular}{ccc}
  \toprule
   \multicolumn{1}{c}{\begin{tabular}[c]{@{}c@{}}Number\\  of nodes\end{tabular}} & \multicolumn{1}{c}{\begin{tabular}[c]{@{}c@{}}Single\\ replication\end{tabular}} & \multicolumn{1}{c}{\begin{tabular}[c]{@{}c@{}}Mini-batch\\ of 10\end{tabular}}\\
   \midrule
    1000 & 0.75s & 1.27s\\
    5000 & 2.89s & 4.14s\\
    10000 & 5.49s & 7.45s\\
    50000 & 26.31s & 34.83s\\
    100000 & 53.12s & 70.27s\\
    500000 & 291.01s & 372.68s\\ 
    \bottomrule
  \end{tabular}
\end{table}

We employ the multiprocessing parallel computing technique and each replication is executed on a separate CPU core. As the communication and data transmission between the parent and children processes takes some time, the run time of multiple replications is higher than that of a  single replication. However, as can be seen in the tables, the run time of a mini-batch of 10 is much lower than 10 times that of a single replication, suggesting that we should take advantage of the multi-core capability in the optimization process.

\subsection{Performance of the Optimization Algorithms}\label{subsec:performance-opt}

Remember that in the Stage 1 of the optimization procedure, we solve a $L_1$ regularized stochastic optimization problem (\ref{eq:L1}). Both the SSGD and FISTA algorithms may be used. We first compare the performance of these two algorithms on the inventory optimization problem with different scales. Notice that the subgradient of the $L_1$ norm (${{{\left\| S \right\|}_1}}$) is needed in the implementation of the SSGD algorithm and we use the sign function (${\mathop{\rm sign}} (S)$). To provide a fair comparison, for both algorithms, we solve the same test problems used in Section \ref{subsec:performance-simu&grad}, we apply the same mini-batch of 10 replications to evaluate the cost and gradient in each epoch, and we use the same regularization parameter ($\lambda$) and the step-size (${t_k})$. Furthermore, because the two algorithms may stop at different epochs, we report the results based on fixed numbers of epochs.

\begin{figure}[ht]
\begin{center}
\includegraphics[scale=0.2]{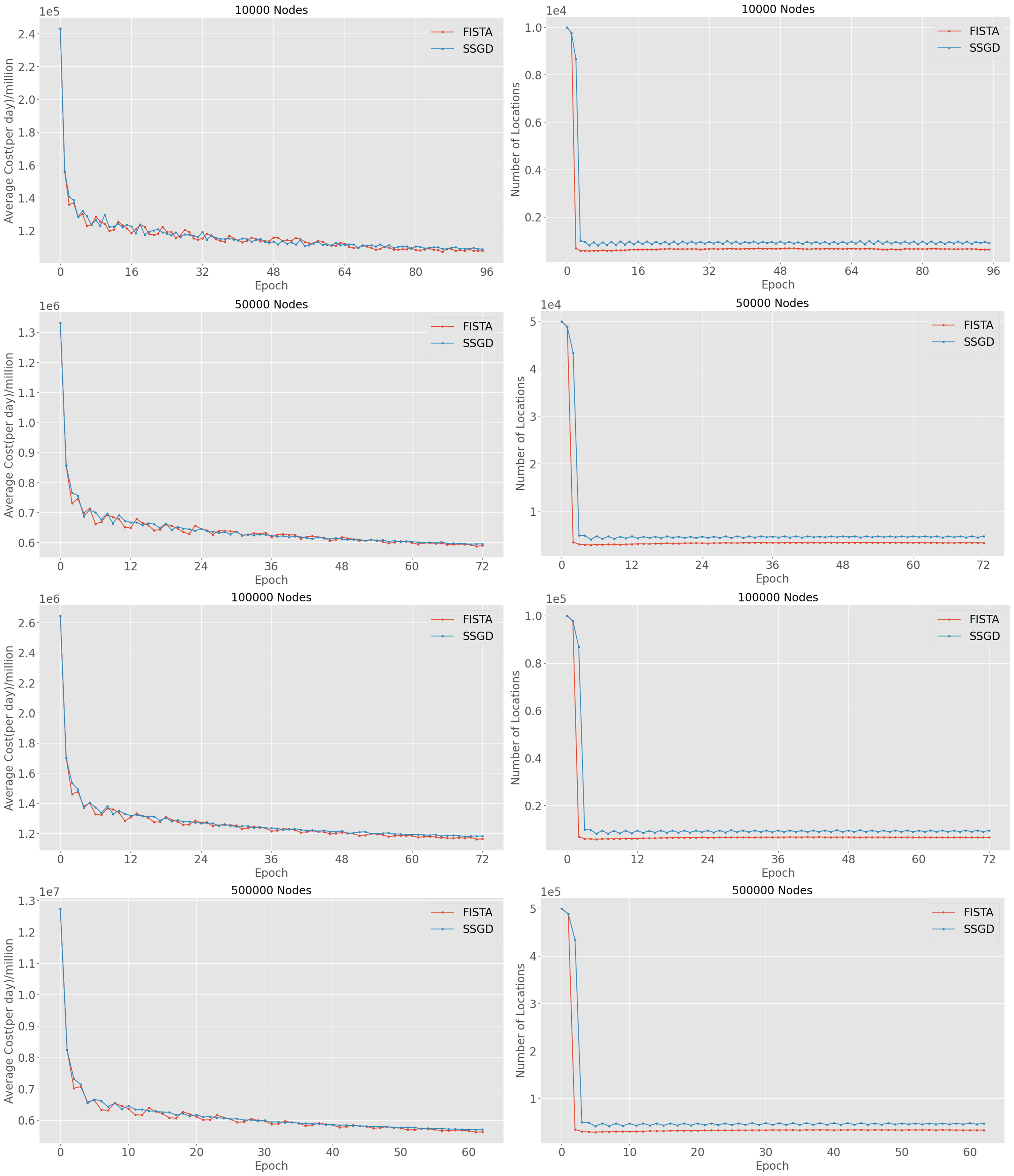}
\caption{Optimization Curve with Fixed Epochs}
\label{fig:Opt-Curve}
\end{center}
\end{figure}

The results are shown in Figure \ref{fig:Opt-Curve}. The four rows represent four different problem scales, i.e., $10,000$, $50,000$, $100,000$ and $500,000$ nodes, and the two columns shows the objective values (i.e., the average total cost) and the sparsity of the solutions (i.e., the numbers of non-zero base-stock levels). From these plots, we see that, in terms of the objective values, both algorithms have similar performance and the FISTA algorithm finds slightly better solutions. In terms of the sparsity of the solutions, the FISTA algorithm tends to find solutions that are more sparse. Therefore, we conclude that in the Stage 1, the FISTA algorithm is a preferred one among the two.

Notice that the solution of Stage 1 provides the information not only on where to keep inventory but also on the base-stock levels. However, in our optimization procedure, we only keep the location information and suggest to use Stage 2 optimization to optimize the base-stock levels. To test the advantage of a second stage, we conduct additional simulation experiments to estimate the objective values of the solutions of Stage 1 and Stage 2 and report the results in  Table \ref{tab:cost improvement}. From the results it is clear that the Stage 2 solution is clearly better than the Stage 1 solution for all four problems, validating the use of the Stage 2 to further improve the quality of the solution from Stage 1.

\begin{table}[ht]\small
  \caption{Cost Improvement in Stage 2}
  \label{tab:cost improvement}
  \centering
  \renewcommand{\arraystretch}{1.1}
  \begin{tabular}{ccccc}
  \toprule
   Number of nodes & 10000 & 50000 & 100000 & 500000\\
   \midrule
   Cost improvement & $3.38\%$ & $2.07\%$ & $3.32\%$ & $2.17\%$ \\
    \bottomrule
  \end{tabular}
\end{table}

In Table \ref{tab:run time} we report the computational times of a typical run of the two-stage optimization procedure, where the FISTA is used in Stage 1 and the SGD is used in Stage 2 and both algorithms are run until the termination conditions are satisfied. From the table, we see that a large-scale inventory optimization problem with 50,000 nodes may be solved in about 2 hours and a very large-scale problem with 500,000 nodes may be solved in about a day. Since inventory decisions are not real-time decisions and they are updated infrequently (for instance, once every six months), these run times are in general acceptable.

\begin{table}[ht]\small
  \caption{Run Time of the Two-Stage Optimization Procedure}
  \label{tab:run time}
  \centering
  \renewcommand{\arraystretch}{1.1}
  \begin{tabular}{ccccc}
  \toprule
   Number of nodes & 10000 & 50000 & 100000 & 500000\\
   \midrule
   Stage 1 & 25.14 min & 92.40 min & 198.10 min & 19.36 hr \\
   \midrule
   Stage 2 & 14.37 min & 26.97 min & 59.41 min & 4.98 hr \\
    \bottomrule
  \end{tabular}
\end{table}

\subsection{Comparison with the Guaranteed Service Model}\label{subsec:comparison}

The simulation optimization algorithm appears to work well in Section \ref{subsec:performance-opt}. However, without a benchmark algorithm to compare with, it is difficult to know the actual performance of the algorithm. To partly solve this problem, we consider the GS model of \cite{graves2000optimizing}. Notice that the GS model requires the BOM to be a spanning tree with bounded demands and deterministic lead times, which are difficult to satisfy for complex production systems. Nevertheless, it provides an optimization algorithm that can be used as a benchmark. In this subsection we consider three different types of test problems. The first one satisfies all requirements of the GS model, the second one has a more complex network structure, and the third one has random lead times. We compare our algorithm with the algorithm of \cite{graves2000optimizing}. From these three types of test problems, we conclude that when the assumptions of the GS model are satisfied, both algorithms produce solutions with similar objectives and solution sparsity. However, when the assumptions are not satisfied, our algorithm performs significantly better than the GS algorithm, as expected.

Throughout this subsection we assume that the outside demands of all items are independent and normally distributed with different means and variances. To implement the GS model, we need an upper bound for each demand and we set it as an upper quantile of the normal distribution as suggested by \cite{graves2000optimizing}. Moreover, the GS model quotes a committed service time to each external demand. We revise the logic of our simulation model slightly to adopt the same logic to ensure a more fair comparison.

\subsubsection{Spanning Tree.}
\begin{figure}[ht]
\begin{center}
\includegraphics[scale=0.35]{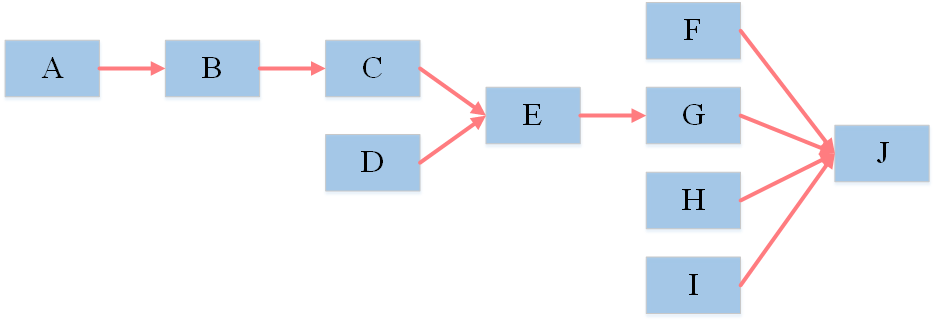}
\caption{Kodak Digital Camera Supply Chain Network}
\label{fig:simple bom1}
\end{center}
\end{figure}

We start with a simple example of Kodak digital camera supply chain pictured in Figure \ref{fig:simple bom1} \citep{snyder2011fundamentals}. The lead times and holding costs are listed in Table \ref{tab:kodak information}. 

\begin{table}[ht]\small
  \caption{Lead Time Information and Holding Cost Coefficients of Kodak Supply Chain}
  \label{tab:kodak information}
  \centering
  \renewcommand{\arraystretch}{1.1}
  \begin{tabular}{m{1cm}<{\centering}m{1cm}<{\centering}m{1cm}<{\centering}m{1cm}<{\centering}m{1cm}<{\centering}m{1cm}<{\centering}m{1cm}<{\centering}m{1cm}<{\centering}m{1cm}<{\centering}m{1cm}<{\centering}m{1cm}<{\centering}}
  \toprule
     & A & B & C & D & E & F & G & H & I & J\\
   \midrule
    ${{l_i}}$ & 2 & 3 & 2 & 4 & 2 & 6 & 3 & 4 & 3 & 2\\
    ${{h_i}}$ & 1 & 3 & 4 & 6 & 12 & 20 & 13 & 8 & 4 & 50\\
    \bottomrule
  \end{tabular}
\end{table}

We apply both algorithms and the optimal base-stock levels found by both algorithms are shown in Table \ref{tab:base-stock1}. As can be seen from the table, both algorithms produce solutions with similar costs and sparsity, with the GS algorithm performs slightly better than our algorithm.

\begin{table}[ht]\small
  \caption{Optimization Results for Kodak Supply Chain Network}
  \label{tab:base-stock1}
  \centering
  \renewcommand{\arraystretch}{1.1}
  \begin{tabular}{m{2.0cm}<{\centering}m{1.1cm}<{\centering}m{1.1cm}<{\centering}m{1.2cm}<{\centering}m{1.2cm}<{\centering}m{1.1cm}<{\centering}m{0.7cm}<{\centering}m{0.8cm}<{\centering}m{0.7cm}<{\centering}m{0.5cm}<{\centering}m{0.5cm}<{\centering}m{2.0cm}<{\centering}}
  \toprule

      \multirow{2}{*}{Algorithm} & \multicolumn{10}{c}{Base-Stock Level}  &  \multirow{2}{*}{Cost}\\
     \cmidrule(l){2-11}
     & A & B & C & D & E & F & G & H & I & J & \\
   \midrule
    GS & 109.70 & 0 & 202.00 & 156.37 & 0 & 0 & 0 & 0 & 0 & 0 & $1.61 \times {10^5}$\\
    Our & 0 & 71.10 & 65.50 & 0 & 102.43 & 0 & 42.47 & 0 & 0 & 0 & $1.65 \times {10^5}$\\
    \bottomrule
  \end{tabular}
\end{table}

To understand the performances of both algorithms for larger-scale problems, we built some larger inventory systems with spanning tree topology. The results of both algorithms are reported in Table \ref{tab:perform-larger}. Again, we see that both algorithms produce solutions with similar costs and sparsity. Furthermore, we notice that, when the number of nodes is large, e.g., 50,000, our algorithm is significantly faster. 
\begin{table}[ht]\small
  \caption{Optimization Performance on Larger Spanning Trees}
  \label{tab:perform-larger}
  \centering

  \renewcommand{\arraystretch}{1.1}
  \begin{tabular}{ccccccc}
  \toprule

    \multirow{2}{*}{\begin{tabular}[c]{@{}c@{}}Number\\  of nodes\end{tabular}}  
       & \multicolumn{3}{c}{GS Model}              & \multicolumn{3}{c}{Our Algorithm}                \\ \cmidrule(l){2-4}\cmidrule(l){5-7}
                         & Cost             & Location  & Time             & Cost                 & Location   & Time           \\
    
   \midrule
    10000 & $3.72 \times {10^{11}}$ & 2637 & 1.92h & $3.97 \times {10^{11}}$ & 2585 & 1.39h\\
    50000 & $1.94 \times {10^{12}}$ & 13215 & 97.70h & $1.87 \times {10^{12}}$ & 15239 & 5.49h\\
    \bottomrule
  \end{tabular}
\end{table}

\subsubsection{Network with Cycles.}

\begin{figure}[ht]
\begin{center}
\includegraphics[scale=0.33]{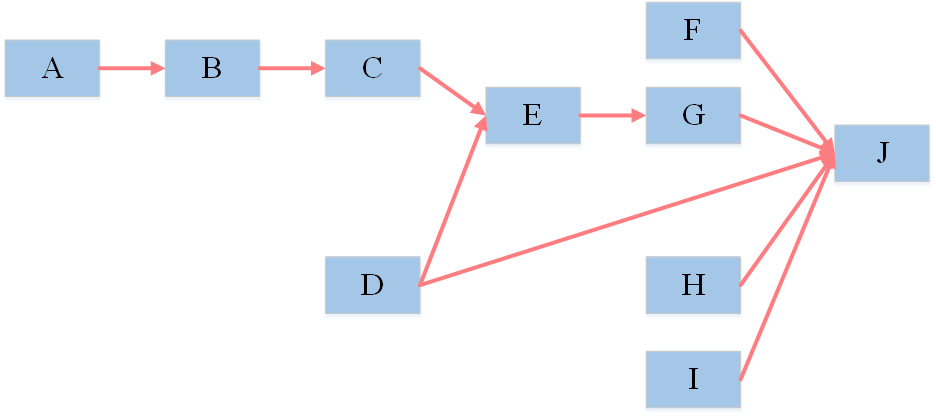}
\caption{Modified Kodak Supply Chain Network}
\label{fig:simple bom2}
\end{center}
\end{figure}

The GS model proposed by \cite{graves2000optimizing} is developed for supply chains that can be modeled as spanning trees. For general systems which may include cycles, this model is no longer applicable. Next, we modify the Kodak network (see Figure \ref{fig:simple bom2}) to explore the influence of the cycles. 

To apply the GS model to the network in Figure \ref{fig:simple bom2}, a simple idea is to cut off an edge and regard it as a spanning tree. Here, we cut off the arc between ``D'' and ``J'' so that the optimal base-stock level is the same as that of Table \ref{tab:base-stock1}. We also obtain the optimal base-stock level by our algorithm, and the results are shown in Table \ref{tab:base-stock2}. From the results we can see that our algorithm finds the solution with two more non-zero base-stock levels and better cost, compared to the solution of the GS model.

\begin{table}[ht]\small
  \caption{Optimization Results for Modified Kodak Supply Chain Network}
  \label{tab:base-stock2}
  \centering
  \renewcommand{\arraystretch}{1.1}
  \begin{tabular}{m{2.0cm}<{\centering}m{1.1cm}<{\centering}m{1.1cm}<{\centering}m{1.2cm}<{\centering}m{1.2cm}<{\centering}m{0.65cm}<{\centering}m{0.5cm}<{\centering}m{0.9cm}<{\centering}m{0.7cm}<{\centering}m{0.5cm}<{\centering}m{0.5cm}<{\centering}m{2.0cm}<{\centering}}
  \toprule
     \multirow{2}{*}{Algorithm} & \multicolumn{10}{c}{Base-Stock Level}  &  \multirow{2}{*}{Cost}\\
     \cmidrule(l){2-11}
     & A & B & C & D & E & F & G & H & I & J & \\
   \midrule
    GS & 109.70 & 0 & 202.00 & 156.37 & 0 & 0 & 0 & 0 & 0 & 0 & $2.27 \times {10^5}$\\
    Our & 58.81 & 63.35 & 55.53 & 0 & 70.02 & 0 & 33.38 & 0 & 0 & 0 & $1.90 \times {10^5}$\\
    \bottomrule
  \end{tabular}
\end{table}

\subsubsection{Random Lead Times.}

The GS model assumes a deterministic lead time for each node in the supply chain. Nevertheless, this is typically not true in practice. For instance, procurement lead times are often quite volatile. We consider a case that the procurement lead times of the raw material nodes in Figure \ref{fig:simple bom1} are random, i.e., ${l_{t,i}} = \max \left( {1,round\left( {N\left( {{\mu _{{l_i}}},0.05{\mu _{{l_i}}}} \right)} \right)} \right)$. For the GS model, we use $\mu_{l_i}$ as the (deterministic) lead time. The optimal solutions found by both algorithms are evaluated based on the simulation model, and the total cost of $100$ periods for the GS model and our algorithm are $1.97 \times {10^6}$ and $9.38 \times {10^5}$, respectively. It is clear that our algorithm performs significantly better than the GS model.

\section{Concluding Remarks}\label{sec:conclud}

Inspired by the similarities between inventory optimization and recurrent neural networks, in this paper we develop a new framework for inventory simulation and inventory optimization for large-scale inventory optimization problems on complex production networks. The framework utilizes modern computational tools to reduce the computational time drastically, so that large-scale problems with tens to hundreds of thousands nodes may be solved in a reasonable amount of time.

There are a few directions that one may further extend our work. First, many practical inventory problems have constraints, e.g., capacity constraints or fill-rate constraints. It is interesting to consider how to incorporate these constraints in the framework. Second, the similarity between RNN and dynamic simulation holds not only for inventory simulation, but also for general periodically reviewed dynamic systems. It is interesting to further study how the framework developed in this paper may be extended to simulation and optimization of general dynamic systems.


%
\begin{APPENDIX}{Derivation of the Inventory Position (Equation \ref{eq.IPti-recursive})}
By definition, the inventory position equals the on-hand inventory plus the on-order quantity minus the backorders (demand that have occurred but have not been satisfied), i.e.,
\begin{equation}
I{P_{t,i}} = I_{t-1,i} + O{O_{t,i}} - \left( {B_{t-1,i}^{out} + B_{t - 1,i}^{in}} \right),
\label{eq.IPti-def}
\end{equation}
where $O{O_{t,i}}$ denotes the on-order quantity of item $i$ at period $t$, ${B_{t - 1,i}^{out}}$ and ${B_{t - 1,i}^{in}}$ denote the backlogged outside and internal demand of item $i$ at period $t-1$, respectively. The on-order quantity refers to the amount that is ordered but not yet received through either procurement or production, i.e.,
\begin{equation}
O{O_{t,i}} = \sum\limits_{u = 1}^{t - 1} {\left( {{O_{u,i}} - {P_{u,i}}} \right)}.
\label{eq.OO}
\end{equation}

The backlogged internal demand can be calculated by
\begin{equation}
B_{t,i}^{in} = \sum\limits_j {{a_{ij}}O_{t,j}^b}.
\label{eq.Btiin}
\end{equation}
Notice that, unlike outside demands where at least one of the on-hand inventory and outside backlog have to be zero (see Equations \ref{eq.Iti-start} and \ref{eq.Btiout}), it is possible to have positive on-hand inventory and positive internal backlog at the same time, because the backlog of the downstream item may be caused by the shortage of some other components.

According to Equation (\ref{eq.Iti-end}), we have
\begin{eqnarray}
I{P_{t,i}} = I_{t - 1,i}^0 - \sum\limits_j {{a_{ij}}{M_{t - 1,j}}}  + O{O_{t,i}} - B_{t - 1,i}^{out} - B_{t - 1,i}^{in}.
\end{eqnarray}
Based on Equations (\ref{eq.Iti-start}) and (\ref{eq.Btiout}), we have \[I_{t - 1,i}^0 - B_{t - 1,i}^{out} = {I_{t - 2,i}} + {P_{t - 1,i}} - B_{t - 2,i}^{out} - D_{t - 1,i}^{out}.\]
Thus,
\begin{eqnarray*}
I{P_{t,i}} &=& {I_{t - 2,i}} + {P_{t - 1,i}} - B_{t - 2,i}^{out} - D_{t - 1,i}^{out} - \sum\limits_j {{a_{ij}}{M_{t - 1,j}}}  + O{O_{t - 1,i}} + {O_{t - 1,i}} - {P_{t - 1,i}} - \sum\limits_j {{a_{ij}}O_{t - 1,j}^b} \\
 &=& {I_{t - 2,i}} + O{O_{t - 1,i}} - \left( {B_{t - 2,i}^{out} + B_{t - 2,i}^{in}} \right) + B_{t - 2,i}^{in} - D_{t - 1,i}^{out} + {O_{t - 1,i}} - \sum\limits_j {{a_{ij}}\left( {{M_{t - 1,j}} + O_{t - 1,j}^b} \right)} \\
 &=& I{P_{t - 1,i}} + {O_{t - 1,i}} - D_{t - 1,i}^{out} + \sum\limits_j {{a_{ij}}O_{t - 2,j}^b}  - \sum\limits_j {{a_{ij}}\left( {{O_{t - 1,j}} + O_{t - 2,j}^b} \right)} \\
 &=& I{P_{t - 1,i}} + {O_{t - 1,i}} - D_{t - 1,i}^{out} - D_{t - 1,i}^{in}.
\end{eqnarray*}
Therefore, Equation (\ref{eq.IPti-recursive}) holds.

\end{APPENDIX}
%
%



\bibliographystyle{informs2014}
\bibliography{main}

\end{document}